\documentclass[12pt]{article}
\usepackage{amsfonts,amssymb,amsmath,theorem}

\let\eps=\varepsilon \let\kappa=\varkappa 
\let\cal=\mathcal \let\bbb=\mathbb \let\goth=\mathfrak  
\newcommand\abs[1]{\left|#1\right|}              
\newcommand\norm[1]{\left\|#1\right\|}           
\newcommand\qed{\ifhmode\unskip\nobreak\fi\quad  
   \ifmmode\square\else\hbox{$\square$}\fi}      
\newcommand\proofskip{\vspace{
   \theorempostskipamount}}                      
\input umsa.fd \input umsb.fd \input ueuf.fd     

\newtheorem{theorem}{Theorem}[section]
\newtheorem{lemma}[theorem]{Lemma}
\newtheorem{prop}[theorem]{Proposition}
\newtheorem{corol}[theorem]{Corollary}
{\theorembodyfont{\normalfont}\newtheorem{rem}[theorem]{Remark}
                        \newtheorem{defn}[theorem]{Definition}
                        \newtheorem{exm}[theorem]{Example}
                        \newtheorem{I}[theorem]{\unskip}}
\newcommand\proof[1]{\noindent\textit{Proof#1}}

\newcommand\Int{\mathop{\mathrm{Int}}}
\renewcommand\circ{\mspace{2mu}{\mathchar"220E}\mspace{2mu}}

\begin{document}

\begin{center}

\Large\bf T-ENTROPY AND VARIATIONAL PRINCIPLE FOR THE SPECTRAL RADIUS OF TRANSFER AND WEIGHTED
SHIFT OPERATORS

\bigskip\bigskip\normalsize\rm

A.\,B.\ ANTONEVICH\footnote{$^{,2\,}$Supported by the Grant of Polish Minister of Science and
Higher Education N201\,382634}

\smallskip
{\it University of Bialystok, Poland\ \ {\rm\&}\ \ Belarus State University, Belarus\par (e-mail:
antonevich@bsu.by)}

\medskip
V.\,I.\ BAKHTIN

\smallskip
{\it Belarus State University, Belarus (e-mail: bakhtin@tut.by)}

\medskip
A.\,V.\ LEBEDEV\footnotemark{}

\smallskip
{\it University of Bialystok, Poland\ \ {\rm\&}\ \ Belarus State University, Belarus\par (e-mail:
lebedev@bsu.by)}

\end{center}

\renewcommand\abstractname{}
\begin{abstract}

The paper deals with the variational principles for evaluation of the spectral radii of transfer
and weighted shift operators associated with a dynamical system. These variational principles have
been the matter of numerous investigations and the principal results have been achieved in the
situation when the dynamical system is either reversible or it is a topological Markov chain. As
the main summands these principles contain the integrals over invariant measures and the
Kolmogorov--Sinai entropy. In the article we derive the  Variational Principle for an
\emph{arbitrary} dynamical system. It gives the  explicit description of the Legendre dual object
to the spectral potential. It is shown that in general this principle contains not the
Kolmogorov--Sinai entropy but a new invariant of entropy type --- the $t$-entropy.
\end{abstract}

\bigbreak\bigskip

\quad\parbox{13.9cm} {\textbf{Keywords:} {\itshape transfer operator, weighted shift, variational
principle, spectral potential, t-entropy}

\medbreak \textbf{2000 MSC:} 37A35; 47B37; 47C15}

\vspace{5mm} \tableofcontents

\vspace{10mm}

\section{Introduction}\label{1..}

The article is devoted to the investigation of dynamical and metric invariants associated with the
spectral radius of transfer  and weighted shift operators.

Let us consider a compact space  $X$, and let  $\alpha\!:X\to X$ be a continuous mapping. This
mapping generates the dynamical system with  discrete time which we will denote by $(X,\alpha)$.

Among the operators  whose spectral analysis is of prime importance in the dynamical systems theory
are

a) the \emph{shift} operators, that is the operators of the form
\begin{equation}\label{1,,1}
 T_\alpha f(x) = f\bigl(\alpha (x)\bigr), \qquad  f\in F(X),
\end{equation}
where $F(X)$ is a certain functional space,

b) \emph{weighted  shift} operators
\begin{equation}\label{1,,2}
 aT_\alpha f(x)=a(x) f\bigl(\alpha(x)\bigr),  \qquad  f\in F(X),
\end{equation}
where $a$ is a fixed function (weight), (operators \eqref{1,,2} are also called \emph{evolution
operators}), and

c) \emph{transfer} operators (associated with the adjoint operators to weighted shift operators
(see Definition \ref{2..1} and Example \ref{7..6})) among which the most popular one is the
classical Perron-Frobenius operator, that is the operator acting in the space $C(X)$ of continuous
functions on $X$ and having the form
\begin{equation}\label{1,,3}
 A f(x)  = \sum_{y\in\alpha^{-1}(x)} {\psi(y)}f(y),
\end{equation}
where  $\psi\in C(X)$ is fixed. This operator is well defined when $\alpha$ is a local
homeomorphism.

Apart from the `pure' dynamical systems theory these operators have  numerous applications in
mathematical physics and in particular in thermodynamics, stochastic processes and information
theory, investigations of zeta functions and Fredholm determinants,   operator algebras theory,
where they serve as an inexhaustible source of important examples and counterexamples so also as
key  constructive elements of the crossed product type algebras, in the theory of solvability of
functional differential equations,  wavelet analysis etc. We refer to books
\cite{Ruelle1,PP,R4,An-book,AnLeb-book1,KS,ABelL1,AbelL2,B} and recent papers
\cite{F,Rg,K,ABLS,Exel,Exel-Ver,ABL-crossed,D} and the bibliography therein.

Spectral properties of weighted shift and transfer operators and especially the formulae and
methods of calculation of their spectral radii are tightly related to the ergodic and entropy
theory of dynamical systems and variational principles of thermodynamic and informational nature.
Let us recall in brief the spectral radius `life story'.

If  $\alpha\!: X \to X$ is a continuous \emph{invertible} mapping and $m$ is an $\alpha $-invariant
measure on $X$ whose support coincides with $X$ (that is any open set has a nonzero measure), and $
a \in C(X)$,  then  in the space $F(X)=L^p(X,m)$, \ $1\le p\le \infty$, as well as in the
space~$C(X)$ the following formulae for the spectral radius of weighted shift operator (\ref{1,,2})
are valid:
\begin{gather}\label{1,,4}
 \ln r(aT_\alpha)=\max_{\mu\in M_\alpha}\int_X\ln|a(x)|\,d\mu,\\[6pt] \label{1,,5}
 \ln r(aT_\alpha)=\max_{\mu\in EM_\alpha}\int_X\ln|a(x)|\,d\mu .
\end{gather}
Here $M_\alpha$ is the set of all  $\alpha$-invariant probability measures on  $X$ and $EM_\alpha$
is the set of $\alpha$-invariant ergodic measures on  $X$.

The statements of this type are called in the dynamical systems theory and  related fields of
analysis the \emph{variational principles} and we will come across a number  of them in  the
article.

The variational principles \eqref{1,,4}, \eqref{1,,5} were stated by Antonevich and for a number of
concrete situations they have  been proved, for example, in \cite{An1, An2} where one can also find
the corresponding description of the set  $M_\alpha$. In the general form (for an \emph{arbitrary}
homeomorphism~$\alpha$) these principles were established by Lebedev \cite{Leb1} and Kitover
\cite{Kit}. The applications of formula  \eqref{1,,5} to the calculation of the spectral radii of
various weighted shift operators are given in \cite{AnLeb,An-book,AnLeb-book1}.

If $\alpha$ is \emph{not invertible} the results related to the corresponding formulae for the
spectral radius can be divided into two classes. The first class contains the results referred to
the case when the weighted shift operators act in the spaces of $C(X)$ or $L^\infty (X,m)$ type. In
this case formulae \eqref{1,,4}, \eqref{1,,5} preserve their form (see, for example, \cite{Serin}).
The second class contains the results referred to the spaces $L^p (X,m)$, \ $1\le p< \infty$. Here
the situation changes drastically. At this point the deep `entropy' and `stochastic' nature of the
spectrum of weighted shift and transfer operators springs out. Namely a complete  description of
this phenomenon is the goal of the paper.

The starting principal results for this second class have been  achieved by Latushkin and Stepin
\cite{Lat,LatStep,LatStep2} under a rather special assumption on the nature of the
mapping~$\alpha$. Namely, in the case when $\alpha$ is a \emph{topological Markov chain} (in
particular $\alpha$ can be an expanding $k$-sheeted cover of a manifold $X$) they proved  the
following formula for the spectral radius of  operator  \eqref{1,,2} in $L^p(X,m)$, \ $1\le
p<\infty$:
\begin{equation}\label{1,,6}
 \ln r(aT_\alpha)=\sup_{\mu\in M_\alpha} \left(\int_X\ln|a(x)|\,d\mu+ \frac{1}{p}\left[
\int_X \ln\rho(x)\, d\mu + h(\mu)\right]\right),
\end{equation}
where  $ h(\mu)$ is the Kolmogorov-Sinai entropy of the measure  $ \mu$ with respect to the
mapping~$\alpha$ and  $\rho$  is a certain continuous nonnegative function defined by this mapping
and such that $\sum_{y\in\alpha^{-1}(x)}\rho(y) \equiv 1$ for any $x\in X$.

Formulae \eqref{1,,4}--\eqref{1,,6} can be considered as  analogues to the known variational
principle that links  entropy and free energy in thermodynamics, and in these formulae the
logarithm of the spectral radius plays the role of free energy. In connection with the problems
considered it is worth mentioning the works by Maslov where similar relations were investigated for
certain evolution differential equations of the form
 \begin{equation*}
 \frac{du}{dt} + Bu =0.
 \end{equation*}
For example, in the paper \cite{Mas}, which is devoted to quantization of thermodynamics, the
coincidence of free energy of thermodynamic system described by this equation and the minimal
eigenvalue of the operator $B$ is established. Note that this minimal eigenvalue is precisely the
\emph{logarithm of the spectral radius} of the operator $e^{-B}$ defining the evolution of the
system.

The  function $\rho$ arising in (\ref{1,,6}) in fact marks a crucial interrelation between weighted
shift and transfer operators in the situation under consideration. Namely, it can be verified that
in this case
 \begin{equation}\label{1,,7}
 \bigl\|(aT_\alpha)^n\bigr\|_{L^p (X,m)}=\bigl\|A^n\bigr\|^{1/p}_{C(X)} \, ,
 \end{equation}

\medskip\noindent
where $A\!: C(X) \to C(X)$ is a transfer (Perron--Frobenius) operator of the form
 \begin{equation}\label{1,,8}
 Af(x)= \sum_{y\in\alpha^{-1}(x)}\vert a \vert^p \rho f\,(y)
 \end{equation}
(cf.\ \eqref{1,,3}).

Equality (\ref{1,,7}) implies the equality
 \begin{equation}\label{1,,9}
 r(aT_\alpha)=r (A)^{1/p}.
 \end{equation}

The explicit calculation of the norm of $aT_\alpha$ by means of formula (\ref{1,,7}) shows, in
particular, that in $L^p$  the norm of the weighted shift operator $aT_\alpha$ with noninvertible
shift is \emph{not} equal to the maximum of the weight $a$ as for invertible shift  but it is equal
to the maximum of the weight \emph{averaged} over inverse images and namely the speed of averaging
(mixing) under the powers of $\alpha$ is reflected in appearance of the entropy type summand in the
right-hand part of the variational principle (\ref{1,,6}).

In connection with formula (\ref{1,,6}) it is reasonable  to recall the variational principle for
the \emph{topological pressure} established by Ruelle \cite{Ruelle} and Walters \cite{Walt}:
 \begin{equation}\label{1,,10}
 P(\alpha, c) = \sup_{\mu\in M_\alpha}  \left(\int_X c(x)d\mu+ h(\mu)\right),
 \end{equation}
where $c\in C(X)$ is a nonnegative function and $P(\alpha, c)$ is the topological pressure defined
by $\alpha$ and $c$. (We would like to stress here that in contrast to (\ref{1,,6}) $\alpha\! :
X\to X$ in (\ref{1,,10}) is an \emph{arbitrary} continuous mapping).

Comparing formulae (\ref{1,,6}), (\ref{1,,9}), and (\ref{1,,10}) we see that in the case when
$\alpha$ is a topological Markov chain the following  relation is valid
 \begin{equation}\label{1,,11}
 \ln r(aT_\alpha)=\frac{1}{p}\ln r (A) = \frac{1}{p}P(\alpha, \ln \psi),
 \end{equation}
where $\psi =\vert a\vert^p\rho$.

The equality in the right-hand part of (\ref{1,,11}), namely,
 \begin{equation}\label{1,,12}
 \ln r (A) = P(\alpha, \ln \psi),
 \end{equation}
has been known, probably, since \cite{Ruelle1}. In fact, the establishment of the relation between
the spectral radius of weighted shift operators and topological pressure was the essence of the
Latushkin--Stepin work. This link along with the observed relation between the spectral radii of
weighted shift and transfer operators serves as a basement for numerous applications of these
operators and also inspires the investigation of their spectral properties in various functional
spaces and, in particular, in the spaces of smooth functions and vector-functions (see, for
example, \cite{CI,BY,CL,GL,BH,GL1}). We have to stress again that all the mentioned sources deal
only with the case when $\alpha $ is a topological Markov chain.

\begin{rem}\label{1..1}
\ 1) In general (that is for an arbitrary  continuous mapping $\alpha$, and even when $\alpha $ is
a local homeomorphism) the equalities $\ln r(aT_\alpha) = \frac{1}{p}P(\alpha, \ln \psi)$ and $\ln
r (B) = P(\alpha, \ln \psi)$ (see (\ref{1,,11})), (\ref{1,,12})) are not true and (\ref{1,,6}) is
\emph{not} a generalization of the variational principle  (\ref{1,,4}). For example, let us
consider an invertible mapping $\alpha$. Then $\rho \equiv 1$. Let us set $a \equiv 1$, thus $\psi
\equiv 1$ and we have
$$
\ln r(T_\alpha)=0,
$$
while
$$
P(\alpha,0)=h(\alpha),
$$
where $h(\alpha)$ is the topological entropy of the dynamical system $(X,\alpha)$, and in general
$h(\alpha)$ could be equal to any nonnegative number.

\medskip

2) A different proof of the Latushkin-Stepin formulae was obtained in  \cite{LebMasl1,LebMasl2} by
means of newly introduced topological invariants that also gave a number of  estimates for the
spectral radius. In addition it was shown in \cite{LebMasl1,LebMasl2} that the variational
principle  (\ref{1,,4}) and the Latushkin-Stepin result are in a way the `extreme points' of the
situations one could come across when dealing with the calculation of the spectral radius
$r(aT_\alpha)$.

\medskip

3) Note that if $p \to \infty$ then formula  (\ref{1,,6})  transforms into (\ref{1,,4})  and this
agrees with the fact (that has been already noted) that in the spaces of $C(X)$ and $L^\infty (X)$
type the variational principle (\ref{1,,4}) preserves its form.
\end{rem}

In the present article we establish the   variational principles for the spectral radii of transfer
and weighted shift operators for an \emph{arbitrary} dynamical system $(X,\alpha)$. It will be
shown that in general these variational principles contain \emph{not} the Kolmogorov-Sinai entropy
$h(\mu)$ as in (\ref{1,,6}) but a \emph{new} dynamical characteristics which we call
{$t$-\emph{entropy}. In the article $t$-entropy is calculated  \emph{explicitly} and its dynamical
and stochastic meanings are uncovered.

The description of $t$-entropy needs the introduction of a deep Legendre transform `philosophy'
into the process of calculation of the spectral radius.  It will be shown that the variational
principles in fact reflect the Legendre duality for the spectral characteristics of  transfer and
weighted shift operators, and $t$-entropy is precisely  the explicit form of the Legendre dual
object to the logarithm of their spectral radii. In particular, the Kolmogorov-Sinai entropy
$h(\mu)$ arising in the right-hand parts of (\ref{1,,6}) and (\ref{1,,10}) is nothing else than
manifestation of the fact that $h(\mu)$ is the Legendre dual object to the logarithm of the
spectral radius of weighted shift operator associated with  the topological Markov chain  and it is
the Legendre dual object to the topological pressure for an arbitrary $\alpha$ respectively.

We would like to emphasize that the main mathematical basement and ideology of the present paper
take their roots in the papers \cite{ABL-2000-1,ABL-2000-2,ABL-2001,ABL-2005,Bakh}. Namely,
particular cases of $t$-entropy and the corresponding variational principles have been considered
in \cite{ABL-2000-1,ABL-2000-2,ABL-2001,ABL-2005} for the situation when the initial transfer
operators are the conditional expectation operators and the definition of $t$-entropy was
introduced there in a \emph{different} way, and in  \cite{Bakh} the variational principle for the
spectral radius of weighted shift operators with positive weights in  $L^1(X,m)$ and the
corresponding Entropy Statistic Theorem (see Section~\ref{9..} below) were  proved. In this paper
we give a complete `final' picture  of the operator-entropy spectral radius phenomena for the
operators in question. Here we present  a general operator algebraic approach establishing the
variational principles as for the spectral radii of \emph{arbitrary}  transfer so also for weighted
shift type operators. The approach developed exploits only the \emph{algebraic} structure of the
objects considered and  does not lean on the spaces of action of  operators (thus essentially
extending the sphere of  applications which therefore do not have to be linked with functional
operators at all).  The results obtained  also unify all the previous structures and we discuss the
interrelation between earlier and presented here definitions in Section \ref{12..}.

 We  note also that by developing the idea of the mentioned  papers Antonevich
and Zajkowski \cite{AZ} proved the convexity of the logarithm of the spectral radii and Legendre
duality for the sums of weighted shift and transfer operators, and in a number of special cases
evaluated explicitly the corresponding `entropies'.

The article is organized as follows. In the starting Sections~\ref{2..} and~\ref{3..} we introduce
the spectral potential (i.\,e. the logarithm of the spectral radius of a family of transfer
operators), examine its elementary properties, and in particular establish its convexity with
respect to weights, and recall the principal properties of the Legendre transform. On this base in
Section~\ref{4..} we introduce the dual entropy as the Legendre dual object to the spectral
potential and derive  the `Legendre' form of the variational principle we are going to investigate.
Our main goal is to obtain the variational principle and the description of the dual entropy in the
\emph{explicit} form \emph{not leaning} on the Legendre transform. This program is implemented in
the further part of the paper. In Section~\ref{5..} we derive a lower estimate for the spectral
potential which in a natural way leads to the main object of the paper --- $t$-entropy (an explicit
dynamical expression for the dual entropy). The lower estimate obtained in this section is in fact
a \emph{sharp} estimate and this is the essence of the Variational Principle which is formulated in
the model situation in Section~\ref{5..} as well. In the next Section \ref{6..} we reformulate  the
Variational Principle for the general \hbox{$C^*$-algebraic} setting, which in particular provides
us with an essential extension of the range of applications, and in Section~\ref{7..} we present a
series of types of transfer operators and $C^*$-dynamical systems naturally arising within the
frames of the $C^*$-algebraic picture chosen. The proof of the Variational Principle needs the
development of an adequate mathematical basement and we provide this in Sections~\ref{8..} and
\ref{9..} where the principal properties of $t$-entropy are examined and the Entropy Statistic
Theorem is proved. The latter theorem gives  statistical estimates of distribution of empirical
measures by means of $t$-entropy and therefore plays for $t$-entropy the role similar to that the
Shannon-McMillan-Breiman  plays for the Kolmogorov-Sinai entropy. In addition it serves as the key
technical instrument in the proof of the Variational Principle itself which is given in
Section~\ref{10..}. Further in Section~\ref{11..} we strengthen the Variational Principle up to the
case of not only positive but also nonnegative coefficients of transfer operators. The Variational
Principle derived tells us in particular that $t$-entropy plays the principal  role in the spectral
analysis of transfer operators and in Section~\ref{12..} we describe  a number of additional
properties of $t$-entropy that have not been exploited in the proof of the Variational Principle.
Along with this we also present a certain alternative definition of $t$-entropy and discuss the
interrelation between the notion of $t$-entropy introduced in this article and that exploited
previously in~\cite{ABL-2000-1,ABL-2000-2,ABL-2001,ABL-2005}. Finally, in Section~\ref{13..} we
examine the interrelation between transfer operators and weighted shift operators acting in
$L^p$-type spaces and prove the corresponding Variational Principle for the spectral radius of
weighted shift operators.

\section[Spectral potential of a transfer operator. Starting properties]{Spectral
 potential of a transfer operator.\\ Starting properties}\label{2..}

Let us consider a Hausdorff compact space $X$. We denote by  $C(X)$ the algebra of continuous
real-valued functions on~$X$ equipped with the uniform norm. Let  $\alpha\!:X\to X$ be a continuous
mapping. This mapping generates the dynamical system with  discrete time which we will denote by
$(X,\alpha)$.

\begin{defn}\label{2..1}
A linear operator  $A\!:C(X)\to C(X)$ is called a \emph{transfer operator} for the dynamical system
$(X,\alpha)$ if

a) $A$  is positive  (that is it maps nonnegative functions to nonnegative) and

b) it satisfies the  \emph{homological identity}
$$
A\bigl(f\circ \alpha\cdot g\bigr) = fAg, \qquad f,g\in C(X)
$$
\end{defn}

A typical example of a transfer operator is given by the classical Perron--Frobenius operator, that
is an operator of the form
\begin{equation}\label{2,,1}
A_\varphi f(x)  = \sum_{y\in\alpha^{-1}(x)} e^{\varphi(y)}f(y),
\end{equation}
where  $\varphi\in C(X)$ is a certain fixed function. This operator is well defined when $\alpha$
is a local homeomorphism. Clearly it is positive and satisfies the homological identity.

Further examples and detailed discussion of transfer operators is presented in Section~\ref{7..}.

Given a transfer operator $A$ we define a family of operators $A_\varphi\!:C(X) \to C(X)$ depending
on the functional parameter $\varphi\in C(X)$ by means of the formula
$$
 A_\varphi f =A(e^\varphi f).
$$
Evidently, all the operators of this family are transfer operators as well. Let us denote
by~$\lambda(\varphi)$ the logarithm of the spectral radius of $A_\varphi$, that is
$$
\lambda(\varphi)  = \lim_{n\to\infty}\frac{1}{n}\ln \norm{A_\varphi^n}.
$$
The positivity of transfer operator implies that
\begin{equation}\label{2,,2}
\lambda(\varphi)  = \lim_{n\to\infty}\frac{1}{n}\ln \norm{A_\varphi^n\mathbf 1},
\end{equation}
where  \textbf{1} is the unit function on  $X$, and $\norm{f}$ denotes the uniform norm of the
function $f\in C(X)$. The functional $\lambda(\varphi)$ is called the \emph{spectral potential} or
the \emph{spectral exponent} of the transfer operator  $A$ (depending on weather we have in mind
dynamical or  spectral associations). In this paper when dealing with the objects associated with
$\lambda(\varphi)$ we are staying on the platform of dynamical (entropy, thermodynamics,
information, stochastics) point of view and therefore throughout the paper $\lambda(\varphi)$ will
be called the \emph{spectral potential}.

Our goal is investigation of $\lambda(\varphi)$.

The next starting proposition   gives   the list of its principal elementary properties.

\begin{prop}\label{2..2}
The spectral potential\/ $\lambda(\varphi)$ is either identically equal to\/  $-\infty$ on the
whole of\/ $C(X)$ or takes only finite values on\/ $C(X)$ and possesses the following properties:

a) \emph{(monotonicity)} \ if $\varphi\le\psi$, then\/ $\lambda(\varphi)\le \lambda(\psi)$;

b) \emph{(additive homogeneity)} \ $\lambda(\varphi +t) =\lambda(\varphi) +t$ for any\/
$t\in\mathbb R$;

c) \emph{(strong  $\alpha$-invariance)} \ $\lambda(\varphi +\psi) =\lambda(\varphi
+\psi\circ\alpha)$ for all\/  $\varphi,\psi\in C(X)$;

d) \emph{(convexity)} \ $\lambda(t\varphi +(1-t)\psi)\le t\lambda(\varphi) +(1-t)\lambda(\psi)$ for
all\/ $\varphi,\psi\in C(X)$ and\/ $t\in [0,1]$;

e) \emph{(Lipschitz property)} \  $\lambda(\varphi) -\lambda(\psi) \le \norm{\varphi-\psi}$ for
all\/ $\varphi,\psi\in C(X)$.

\noindent
In particular, the spectral potential is continuous.
\end{prop}

\proof.
a)  The monotonicity of $\lambda(\varphi)$ follows from the positivity of $A$.

b)  Let us introduce the notation
\begin{equation*}
 S_n\varphi :=\varphi+\varphi\circ\alpha+\,\cdots\,+\varphi\circ\alpha^{n-1}, \qquad \varphi\in
 C(X).
\end{equation*}
Applying $n$ times the homological identity to the operator  $A_\varphi^n =(Ae^{\varphi})^n$ we
obtain
\begin{equation}\label{2,,3}
A_\varphi^nf =A(e^\varphi A(e^\varphi\dotsm A(e^\varphi f)...)) =A^n\bigl(e^{S_n\varphi}f\bigr).
\end{equation}
By substituting this equality into the definition of the spectral potential \eqref{2,,2} we deduce
the additive homogeneity of $\lambda(\varphi)$.

c)  Note that
\begin{equation*}
S_n(\varphi+\psi\circ\alpha) =S_n(\varphi+\psi) +\psi\circ\alpha^n -\psi.
\end{equation*}
Therefore
\begin{equation*}
e^{-t}(A_{\varphi+\psi})^n\le (A_{\varphi+\psi\circ\alpha})^n \le e^{t}(A_{\varphi+\psi})^n, \qquad
\text{where} \quad t=2\norm\psi.
\end{equation*}
Substituting these inequalities into  \eqref{2,,2} we obtain the strong  $\alpha$-invariance of
$\lambda(\varphi)$.

d) Let us fix a point  $x\in X$ and consider the linear functional  $\mu(f)\! :=A^n\!f(x)$
on~$C(X)$. This functional is positive and therefore by the Riesz theorem it can be identified with
a certain Borel measure on   $X$. By applying the  H\"older inequality (with $p=\frac{1}{t}$ and
$q= \frac{1}{1-t}$) to this measure we obtain
\begin{equation*}
\mu\Bigl(e^{S_n(t\varphi +(1-t)\psi)}\Bigr)\le \Bigl(\mu\bigl(e^{S_n\varphi}\bigr)\Bigr)^{t}
\Bigl(\mu\bigl(e^{S_n\psi}\bigr)\Bigr)^{1-t}, \qquad t\in (0,1).
\end{equation*}
This inequality implies in turn that the functional  $\ln\mu\bigl(e^{S_n\varphi}\bigr) =\ln
[A_\varphi^n\mathbf{1}](x)$ is convex with respect to  $\varphi\in C(X)$. Therefore the spectral
potential
\begin{equation*}
\lambda(\varphi) =\lim_{n\to\infty}\frac{1}{n}\ln\norm{A_\varphi^n\mathbf{1}} =\lim_{n\to
\infty}\frac{1}{n}\sup_{x\in X}\ln[A_\varphi^n\mathbf{1}](x)
\end{equation*}
is convex with respect to  $\varphi$ as well.

e)  The monotonicity and additive homogeneity of the spectral potential imply its Lipschitz
property. Indeed,
\begin{equation*}
\lambda(\varphi) -\lambda(\psi) \le \lambda(\psi+\norm{\varphi-\psi}) -\lambda(\psi) =
\norm{\varphi-\psi}.
\end{equation*}

Finally, the Lipschitz property implies that $\lambda (\varphi)$ is either identically equal to
$-\infty$ on the whole of \ $C(X)$ or takes only finite values on \ $C(X)$. \qed\proofskip

\section{Convex functionals and the Legendre transform}\label{3..}

Proposition \ref{2..2} shows in particular that  the spectral potential is a convex functional
on~$C(X)$. As is known among the standard instruments of investigation of convex functionals is the
Legendre transform. In this section we recall the principal notions and facts related to this
transform (in essence they are borrowed from~\cite{??}).

Let  $f$ be a functional on a real Banach space  $L$ with the values in the extended real straight
line $\bar{\mathbb R} =[-\infty,+\infty]$. The set  $D(f) =\{\varphi\in L\mid f(\varphi)
<+\infty\}$ is called the  \emph{effective domain} of the functional  $f$. The functional  $f$ is
called  \emph{convex}, if for all $\varphi,\psi\in D(f)$ and $t\in [0,1]$ the following inequality
holds
\begin{equation*}
f\bigl(t\varphi +(1-t)\psi\bigr) \le tf(\varphi) +(1-t) f(\psi).
\end{equation*}

The functional  $f$ is called  \emph{lower semicontinuous} if the set  $\{\varphi\in L\mid
f(\varphi) > c\}$ is open for any  $c\in \bar{\mathbb R}$. One can speak  about lower
semicontinuity with respect to the norm topology or with respect to the weak topology on  $L$, but
for convex functionals these properties are equivalent.

Let  $L^*$ be the dual space to  $L$. The functional $f^*\!:L^*\to \bar{\mathbb R}$ that  is
defined on the dual space by the equality
\begin{equation*}
f^*(\mu) =\sup_{\varphi\in L}\bigl(\mu(\varphi) -f(\varphi)\bigr) = \sup_{\varphi\in
D(f)}\bigl(\mu(\varphi) -f(\varphi)\bigr), \qquad \mu\in L^*,
\end{equation*}
is called the  \emph{Legendre dual} to the functional $f$ (or the Legendre transform of $f$). For a
functional $g$ on the dual space the Legendre transform is defined as the functional on the initial
space given by the similar formula:
\begin{equation*}
g^*(\varphi) =\sup_{\mu\in {L}^*}\bigl(\mu(\varphi) -g(\mu)\bigr),\qquad \varphi\in L.
\end{equation*}

\begin{prop}\label{3..1}
Let a functional\/ $f\!:L\to (-\infty,+\infty]$ be not identically equal to\/~$+\infty$. Then

a) the dual functional\/ $f^*$ is convex and lower semicontinuous with respect to  *-weak topology
on the dual space;

b) if the functional\/ $f$ is convex and lower semicontinuous then\/ $f=(f^*)^*$ (the Legendre
transform is involutory);

c) in general\/ $(f^*)^*$ is the maximal convex lower semicontinuous functional that does not
exceed\/ $f$ (the convex hull of\/ $f$).

The analogous statements are valid for functionals\/ $g\!:L^*\to (-\infty,+\infty]$.
\end{prop}

We have already proved that the spectral potential  $\lambda(\varphi)$ is convex and continuous
(see Proposition \ref{2..2}). Therefore it can be represented as the Legendre transform of its
Legendre dual on the dual space  $C^*(X)$. However we will slightly modify the form of the record
of this duality. The matter is that in thermodynamics, information theory and ergodic theory there
is a tradition  to change the sign of the dual to $\lambda(\varphi)$ and the result obtained is
called the  \emph{entropy}.  Following this tradition and, in particular,  the results
of~\cite{ABLS}, we define the \emph{dual entropy} $S(\mu)$ of the spectral potential
$\lambda(\varphi)$ by means of the formula
\begin{equation}\label{3,,1}
S(\mu) \,:=\inf_{\varphi\in C(X)}\bigl(\lambda(\varphi) - \mu(\varphi)\bigr),\qquad \mu\in C^*(X).
\end{equation}

\smallskip\noindent
Since the dual entropy differs from the dual to $\lambda(\varphi)$ functional only by sign it
follows that $S(\mu)$ is concave and upper semicontinuous (with respect to $^*$-weak topology). As
the Legendre transform is involutory  the next equality holds true
\begin{equation}\label{3,,2}
\lambda(\varphi) =\sup_{\mu\in C^*(X)}\bigl(\mu(\varphi) + S(\mu)\bigr).
\end{equation}
This equality is in fact \emph{the simplest form} of the subject of the article --- the
\emph{Variational Principle} for $\lambda(\varphi)$. The foregoing discussion implies the existence
of the dual entropy and its explicit construction by means of formula (\ref{3,,1}) is possible
provided that the spectral potential is known. However, the spectral potential itself is a rather
hard object to investigate. Our principal goal is the \emph{independent} derivation of the
\emph{explicit} formula for the dual entropy, \emph{not leaning}  on the spectral potential. This
formula allows, in particular, to impart a more effective character to the Variational Principle.
The principal result for a model example will be formulated  in Section \ref{5..}, its general
$C^*$-formulation is given in Section \ref{6..} and its complete proof will be given in Section
\ref{10..}.

\section{Dynamical potentials and dual entropy}\label{4..}

Certain useful information on the dual entropy  $S(\mu)$ can be extracted by means of the Legendre
duality from the properties of the spectral potential already proved. In essence this job was
implemented in   \cite{Bakh1}. For completeness of description we present a part of the material of
this article in this section.

Let, as above, $X$ be a Hausdorff compact space and $\alpha\!:X\to X$ be a continuous mapping. We
denote by  $M\subset C^*(X)$ the set of all linear positive normalized functionals on~$C(X)$ (that
is linear functionals that take nonnegative values on nonnegative functions and are equal to $1$ on
the unit function). By the Riesz theorem these functionals are bijectively identified with regular
probability Borel measures on~$X$, and the elements of~$C^*(X)$ are identified with regular
real-valued Borel measures on~$X$. Therefore with a slight abuse of language all the elements of
$C^*(X)$ will be referred to as measures.

A measure  $\mu\in M$ is called  $\alpha$-invariant if  $\mu(f\circ\alpha)  = \mu(f)$ for all
functions $f\in C(X)$. This is equivalent to the identity $\mu(\alpha^{-1}(G)) = \mu(G)$ for all
Borel subsets  $G\subset X$. The collection of all  $\alpha$-invariant measures from  $M$ will be
denoted by  $M_\alpha$.

Recall that according to  \cite{Bakh1} a \emph{dynamical potential} of a dynamical system
$(X,\alpha)$ is an \emph{arbitrary} real-valued functional  $\lambda(\varphi)$ on $C(X)$,
possessing the following  four properties

a) \emph{monotonicity} (if  $\varphi\le\psi$, then  $\lambda(\varphi)\le \lambda(\psi)$);

b) \emph{additive homogeneity} ($\lambda(\varphi +t) =\lambda(\varphi) +t$ for any  $t\in\mathbb
R$);

c) \emph{strong\/ $\alpha$-invariance} ($\lambda(\varphi +\psi) =\lambda(\varphi +\psi\circ\alpha)$
for all  $\varphi,\psi\in C(X)$);

d)  \emph{convexity} with respect to   $\varphi$.

Proposition  \ref{2..2} means precisely that the spectral potential of a transfer operator is a
particular case  (and a typical example) of a dynamical potential.

Let $\lambda(\varphi)$ be \emph{any} dynamical potential on  $C(X)$. We define the dual
entropy~$S(\mu)$ for this dynamical potential just as for the spectral potential by means of the
same formula~\eqref{3,,1}. Clearly formula~\eqref{3,,2} holds true as well.

An \emph{equilibrium measure}, corresponding to a function  $\varphi\in C(X)$, is an arbitrary
subgradient of the functional  $\lambda(\varphi)$ at the point  $\varphi$ (in other words it is a
linear functional $\mu\!:C(X) \to\bbb R$ such that $\lambda(\varphi +\psi)- \lambda(\varphi)\ge
\mu(\psi)$ for all  $\psi\in C(X)$). Evidently the set of all equilibrium measures corresponding to
a certain function $\varphi$ is convex and closed  (with respect to the  $^*$-weak topology). This
set is nonempty by the  convex analysis theorem on the existence of a supporting hyperplane. It
consists of a unique measure~$\mu$ if and only if there exists the G\^ateaux derivative~$\lambda'
(\varphi)$. In this case $\mu =\lambda' (\varphi)$.

The definitions of the dual entropy and equilibrium measure imply

\begin{prop}\label{4..1}
For any function\/ $\varphi\in C(X)$ and any measure\/ $\mu\in C^*(X)$ the Young inequality holds
true
 \begin{equation} \label{4,,1}
 S(\mu)\le \lambda(\varphi) -\mu(\varphi).
 \end{equation}
This inequality turns into equality iff\/ $\mu$ is an equilibrium measure corresponding
to\/~$\varphi$.
\end{prop}

\begin{prop}\label{4..2}
The effective domain of\/ $S$ is contained in\/ $M_\alpha$, that is if\/ $S(\mu)>-\infty$, then\/
$\mu$ is a probability and\/ $\alpha$-invariant measure. In particular this is true for all
equilibrium measures.
\end{prop}

\proof{} is implemented by contradiction. By the additive homogeneity of the dynamical potential
for all $t\in\mathbb R$ we have
$$
 S(\mu)\le \lambda(t) -\mu(t) =\lambda(0)+t(1-\mu(1)).
$$
Therefore, if  $\mu(1)\ne 1$ then $S(\mu) =-\infty$.

Suppose that  $\mu(\varphi)<0$ for some nonnegative function  $\varphi\in C(X)$. By the
monotonicity of the dynamical potential  for all  $t>0$ the following inequalities hold
$$
 S(\mu)\le\lambda (-t\varphi) -\mu(-t\varphi)\le \lambda(0)+ t\mu(\varphi),
$$
and again  $S(\mu) =-\infty$.

The invariance of a measure $\mu$ is equivalent to the identity $\mu(\varphi) \equiv
\mu(\varphi\circ \alpha)$,\ \ $\varphi \in C(X)$. Suppose that
$\mu(\varphi)\ne\mu(\varphi\circ\alpha)$ for a certain function $\varphi$. By the strong invariance
of the dynamical potential we have
$$
\lambda(t\varphi-t\varphi\circ\alpha) \equiv \lambda(0).
$$
Therefore,
$$
S(\mu)\le\lambda(t\varphi-t\varphi \circ \alpha) -\mu(t\varphi -t\varphi\circ\alpha) =\lambda(0)-t
(\mu(\varphi) - \mu(\varphi \circ\alpha))$$ and thus  $S(\mu) =-\infty$.

Finally, if  $\mu$ is an equilibrium measure corresponding to a function $\varphi\in C(X)$ then the
foregoing proposition implies $S(\mu) =\lambda(\varphi) -\mu(\varphi) >-\infty$. \qed \proofskip

The proposition just proved shows that it suffice to define the dual functional (dual entropy) only
on invariant probability measures. This will be done in the next section.

\begin{prop}\label{4..3}
If\/ $S(\mu) >-\infty$, then\/ $\mu$ belongs to the closure of the set of equilibrium measures\/
\textup{(}with respect to the norm of the space\/ $C^*(X)$\textup{)}.
\end{prop}

\proof. A particular variant of the Bishop--Phelps theorem \cite{Bishop} states that if
$\lambda(\varphi)$ is a continuous convex functional on a Banach space $L$, \,$\mu\in L^*$, and the
difference  $\lambda -\mu$ is bounded from below then the distance from  $\mu$ to the set of
subgradients of $\lambda$ is zero. \qed \proofskip

Proposition~\ref{4..2} implies that the supremum in \eqref{3,,2} is attained on the set of
invariant probability measures  $M_\alpha$. Therefore, every dynamical potential has the form
\begin{equation}\label{4,,2}
 \lambda(\varphi) =\sup_{\mu\in M_\alpha} \bigl(\mu(\varphi) +S(\mu)\bigr),
\end{equation}
where  $S(\mu)$ is the corresponding dual entropy.

Proposition~\ref{4..1} implies in turn that supremum in~\eqref{4,,2} is in fact maximum that is
\begin{equation}\label{4,,3}
\lambda(\varphi) =\max_{\mu\in M_\alpha} \bigl(\mu(\varphi) +S(\mu)\bigr),
\end{equation}
and this maximum is attained precisely on equilibrium measures, corresponding to the function
$\varphi$.

The forgoing observation implies in addition that the uniqueness  of an extremal measure
in~\eqref{4,,3} is equivalent to the existence of the G\^ateaux derivative $\lambda'(\varphi)$.

Let us present one more statement describing the structure of dynamical potentials.

\begin{prop}\label{4..4}
If\/ $S(\mu)$ is an arbitrary function which is  bounded from above on\/ $M_\alpha$, then
formula\/~\eqref{4,,2} defines a dynamical potential.
\end{prop}

\proof{} of this fact reduces to a trivial verification of all the conditions of the definition of
the dynamical potential.  \qed \proofskip

\begin{rem}
Since the dual entropy satisfies the inequality  $S(\mu)\le \lambda(0)$, it follows that formula
\eqref{4,,2} implies the estimate
$$
  \lambda(\varphi) \le \sup_{\mu\in M_\alpha}\mu(\varphi) +\lambda(0).
$$
On the other hand, if we take  $S(\mu) \equiv \lambda(0)$, then~\eqref{4,,2} defines the dynamical
potential
$$
\lambda_0(\varphi) :=\sup_{\mu\in M_\alpha} \mu(\varphi) +\lambda(0).
$$
The previous estimate  shows that this dynamical potential is maximal among all the dynamical
potentials with the fixed value~$\lambda(0)$.
\end{rem}

\section[Lower estimate for the spectral potential, $\boldsymbol{t}$-entropy,\\ and formulation of
the Variational Principle in the model setting]{Lower estimate for the spectral potential,\\
$\boldsymbol{t}$-entropy, and formulation of the Variational\\ Principle in the model
setting}\label{5..}

In the previous section we  defined the dual entropy  of an arbitrary dynamical potential by means
of the Legendre transform. Let us consider now  a particular case when the dual entropy is
generated by the spectral potential of a transfer operator. Our goal is to derive a direct
definition for this dual entropy not leaning on the Legendre duality but only on the properties of
the initial dynamical system and the transfer operator chosen. As it was already noted Proposition
\ref{4..2} implies that it suffice to define the dual entropy only on invariant probability
measures.

Let $A\!:C(X)\to C(X)$ be a transfer operator with the spectral potential $\lambda(\varphi)$,
associated with a dynamical system $(X,\alpha)$, where $\alpha$ is a continuous mapping of a
compact space $X$.

\medskip

\noindent
{\bf Lower estimate. Preliminaries.\,}  We start with a certain heuristic reasoning leading us to
the definition we are looking for. Let us show how for any invariant probability measure $\mu$ one
can find a number $C(\mu)$ satisfying inequality (\ref{4,,1}), that is
$$
C(\mu) \le \lambda(\varphi) -\mu(\varphi).
$$

Recall the notation
\begin{equation}\label{5,,1}
 S_n\varphi:=\varphi+\varphi\circ\alpha+\dotsm+\varphi\circ\alpha^{n-1}.
\end{equation}
Consider the expression
$$
 \ln \bigl\|A^n_\varphi\mathbf{1}\bigr\| = \max _{x\in X} \,\ln\left[A^n(e^{S_n(\varphi)})\right](x).
$$

Given an invariant probability measure  $\mu$  we have to find a number $C(\mu)$ such that for
every $\varphi$ the next lower estimate holds
$$
 \frac{1}{n}\ln \bigl\| A^n_\varphi\mathbf{1}\bigr\|\geq \int_X \varphi\, d\mu + C(\mu).
$$

By the invariance of $\mu$ we have
\begin{equation} \label{5,,2}
 \int_X S_n\varphi\,d\mu = n\int_X \varphi\, d\mu.
\end{equation}
Let us replace the integral in the left-hand part of (\ref{5,,2}) by an integral sum (with an
accuracy up to $\eps$). To implement this we take a (sufficiently fine) partition of $X$
$$
D = \{ D_1, \ldots, D_k\}, \qquad X = \coprod_{i=1}^k D_i,
$$
where each $D_i\subset X$ is a Borel set such that the oscillation of the function $S_n\varphi$ on
$D_i$ does not exceed $\eps$. Set
$$
S_n\varphi(D_i) := \sup\{\,S_n\varphi(x)\mid x \in  D_i\,\}.
$$
Then
$$
\sum_{i=1}^k S_n\varphi(D_i)\mu (D_i) -\eps \leq \int_X S_n\varphi\,d\mu \leq \sum_{i=1}^k
S_n\varphi(D_i)\mu (D_i).
$$

Let us estimate the value of $\ln [A^n_\varphi\mathbf{1}](x)$, where $x\in X$ is arbitrary. Denote
by $g_i$ the index function of the set $D_i$. Then the following inequalities hold true
\begin{gather*}
\ln \bigl[A_\varphi^n\mathbf 1\bigr](x) =
\ln\sum_{i=1}^k\bigl[A^n\bigl(e^{S_n\varphi}g_i\bigr)\bigr](x)  \\[3pt]
\geq \ln\sum_{i=1}^ke^{S_n\varphi(D_i)- \eps}[A^n g_i](x)=
 \ln\left[\sum_{i=1}^ke^{S_n\varphi(D_i)}[A^n g_i](x_0)\right]- \eps.
  \end{gather*}
To estimate the logarithm  first we estimate the sum obtained
$$
 \sum_{i=1}^ke^{S_n\varphi(D_i)}\bigl[A^n g_i\bigr](x) \ge\sum_{i=1}^k  \mu(D_i)
 e^{S_n\varphi(D_i)} \frac{[A^n g_i](x)}{\mu (D_i)},
$$
where we assume that if $\mu(D_i)  = 0$ for a certain $D_i$, then we set the corresponding summand
in the right-hand part to be zero (regardless the value of~$[A^ng_i](x)$). Now we exploit the
concavity of the logarithm function and \eqref{5,,2}:
\begin{gather*}
 \ln \bigl[A_\varphi^n\mathbf 1\bigr](x) \geq
 \sum_{\mu(D_i)\ne 0}\mu(D_i)\ln\frac{e^{S_n\varphi(D_i)}[A^n g_i](x)}{\mu(D_i)} -\eps \displaybreak[0]\\[3pt]
 =\sum_{\mu(D_i)\ne 0}\mu(D_i)S_n\varphi(D_i) + \sum_{\mu(D_i)\ne 0}\mu(D_i)\ln\frac{[A^n
 g_i](x_0)}{\mu(D_i)} -\eps \displaybreak[0] \\[3pt]
 \geq n  \int_X \varphi\, d\mu + \sum_{\mu(D_i)\ne 0}\mu(D_i)\ln\frac{[A^n g_i](x)}{\mu(D_i)}- \eps.
\end{gather*}
This implies the inequality
\begin{equation} \label{5,,3}
 \frac{1}{n}\ln\bigl\|A_\varphi^n\mathbf 1\bigr\| \geq \int_X \varphi\,d\mu +
 \frac{1}{n}\sup_{x\in X} \sum_{\mu(D_i)\ne 0}\mu(D_i)\ln\frac{[A^ng_i](x)}{\mu(D_i)}-\frac{\eps}{n}.
\end{equation}

The essential point here is that the second summand in the right-hand part of the inequality does
not depend on $\varphi$. Since we are seeking the limit inequality as $n\to \infty$ we naturally
arrive at the introduction of the new value:
\begin{equation}\label{5,,4}
 C(\mu) :=\inf_{n\in\mathbb N} \inf_{D}\frac{1}{n}\sup_{x\in X}
\sum_{\mu(D_i)\ne 0}\mu(D_i)\ln\frac{[A^ng_i](x)}{\mu(D_i)}.
\end{equation}
Now tending $n$  in (\ref{5,,3}) to infinity and using (\ref{5,,4}) and the arbitrariness of
$\eps>0$ we obtain the estimate of the form required
\begin{equation}\label{5,,5}
 \lambda(\varphi)\ge  \int_X \varphi\, d\mu +C(\mu).
\end{equation}

Let us make two observations regarding  the reasoning presented.

\medskip

1) The index functions $ g_i$ exploited in the calculation are \emph{discontinuous}. As a rule, in
examples the action of the operator $A$ is naturally defined on these functions as well, but in
general  the expression $ A^n g_i$ is \emph{not} defined an therefore the reasoning presented needs
a modification. Because of this in what follows we will consider certain partitions of unity in the
algebra $C(X)$ instead of partitions of the space $X$.

\medskip

2) Expression (\ref{5,,4}) for the constant that satisfies  estimate (\ref{5,,5}) is \emph{not} the
best one. One can improve it by means of the following consideration. The functional $\delta_{x}(f)
= f(x)$ can be identified with the probability measure concentrated at the point $x$ and so
estimate (\ref{5,,3}) can be considered as an estimate obtained by means of the measures
$\delta_x$. It will be shown below that one can replace $[A^ng_i](x) =\delta_x(A^ng_i)$ in
(\ref{5,,3}) by $m(A^ng_i)$, where $m$ is an arbitrary probability measure on $X$. Therefore the
supremum over $x\in X$ in (\ref{5,,3}) can be replaced by the analogous supremum over all
probability measures on $X$, which brings us to a definition of a certain new value $\tau(\mu)$
greater than $C(\mu)$. The principal result of this article is the proof of the fact that the value
of the constant $\tau(\mu)$ obtained in this manner is the maximal possible one, which in turn
implies the Variational Principle.

From the geometrical point of view the situation here is the following. The set $M$ of all
probability measures on $X$ is convex and the measures $\delta_{x}$ are the extreme pints of this
set. The expression obtained in the corresponding estimate depends concavely on measure  $m \in M$.
Therefore the supremum over the set $M$ in the general case is greater than the supremum over the
extreme points of this set.

\medskip

Now we pass to the strict reasoning.

\medskip

By a \emph{partition of unity} in the algebra  $C(X)$ we mean any finite set $D =
\{g_1,\dotsс,g_k\}$ consisting of nonnegative functions $g_i\in C(X)$ satisfying the identity $g_1+
\dots +g_k \equiv 1$.

Our first step is the next technical lemma which gives a key to  derivation  of a lower estimate
for the spectral potential.

\begin{lemma}\label{5..1}
For any $n\in \mathbb{N}$ and $ \mu\in M_\alpha$  the following estimate holds true
\begin{equation}\label{5,,6}
 \ln \norm{A_\varphi^n \mathbf 1}\ge n\mu(\varphi) +\inf_D\sup_{m\in M} \sum_{g\in D}
\mu(g) \ln\frac{m(A^ng)}{\mu(g)}
\end{equation}
{\em (if $\mu(g) =0$ then we put the corresponding summand in \eqref{5,,6} to be zero)}.
\end{lemma}

\proof. Take any positive integer $n$ and a number  $\eps>0$. For these numbers let us choose a
(fine) partition of unity~$D$ in the algebra~$C(X)$ such that the oscillation of the function
$S_n\varphi =\varphi+\varphi\circ\alpha+\dotsm+\varphi\circ\alpha^{n-1}$ on the support of each
function $g\in D$ does not exceed~$\eps$. Set
$$
S_n\varphi(g) := \sup\{\,S_n\varphi(x)\mid g(x)\ne 0\,\}.
$$

\medskip\noindent
Equality  \eqref{2,,3} and concavity of the logarithm function imply the following inequalities for
all probability measures  $\mu\in M_\alpha$ and $m\in M$:
\begin{gather*}
 \eps+\ln m(A_\varphi^n\mathbf 1)
  = \eps +\ln\sum_{g\in D}m\bigl(A^n\bigl(e^{S_n\varphi}g\bigr)\bigr)
 \ge \ln\sum_{g\in D}e^{S_n\varphi(g)}m(A^ng)\displaybreak[0]\\[3pt]
 \ge\ln\sum_{\mu(g)\ne 0}\mu(g)\frac{e^{S_n\varphi(g)}m(A^ng)}{\mu(g)}
 \ge \sum_{\mu(g)\ne 0}\mu(g)\ln\frac{e^{S_n\varphi(g)}m(A^ng)}{\mu(g)}\displaybreak[0]\\[3pt]
 =\sum_{\mu(g)\ne 0}\mu(gS_n\varphi (g)) + \sum_{g\in D}\mu(g)\ln\frac{m(A^ng)}{\mu(g)}
  \ge \sum_{\mu(g)\ne 0}\mu(gS_n\varphi)
 + \sum_{\mu(g)\ne 0}\mu(g)\ln\frac{m(A^ng)}{\mu(g)}\\[3pt]
 =\mu(S_n\varphi) + \sum_{\mu(g)\ne 0}\mu(g)\ln\frac{m(A^ng)}{\mu(g)}.
\end{gather*}
Passing in these inequalities
 to the supremum over  $m\in M$ we obtain the inequality
\begin{equation*}
\eps+ \ln \norm{A_\varphi^n \mathbf 1}\ge \mu(S_n\varphi) +\sup_{m\in M} \sum_{\mu(g)\ne 0} \mu(g)
\ln\frac{m(A^ng)}{\mu(g)}.
\end{equation*}
This inequality along with $\alpha$-invariance of $\mu$ and the arbitrariness of $\eps>0$ implies
\begin{equation*}
\ln \norm{A_\varphi^n \mathbf 1}\ge n\mu(\varphi) +\inf_D\sup_{m\in M} \sum_{\mu(g)\ne 0} \mu(g)
\ln\frac{m(A^ng)}{\mu(g)}. \qed
\end{equation*}

\begin{corol}\label{5..2}
For any $\mu\in M_\alpha$ we have
\begin{equation}\label{5,,7}
 \lambda(\varphi)\ge \mu(\varphi)
 +\inf_{n\in\mathbb N}\frac{1}{n}\inf_D\sup_{m\in M} \sum_{g\in D} \mu(g) \ln\frac{m(A^ng)}{\mu(g)}
\end{equation}
\end{corol}

\proof. Divide \eqref{5,,6} by $n$ and tend $n$ to infinity. \qed \proofskip

The lower estimate for the spectral potential obtained in Corollary \ref{5..2} makes it natural the
introduction of the following object ($t$-entropy), which in essence describes the expression in
the right-hand part of  \eqref{5,,7}.

\begin{defn}\label{5..3}
\emph{$T$-entropy} is the  functional $\tau$  on $M_\alpha$ such that its value at $\mu\in
M_\alpha$ is defined by the following formulae
\begin{gather}\label{5,,8}
 \tau(\mu)  := \inf_{n\in\mathbb N}\frac{\tau_n(\mu)}{n}\,,\qquad
 \tau_n(\mu)  := \inf_D\tau_n(\mu,D),\\[6pt]
 \tau_n(\mu,D)  := \sup_{m\in M}\sum_{g\in D}\mu(g) \ln\frac{m(A^ng)}{\mu(g)}\,. \label{5,,9}
\end{gather}
The infimum in \eqref{5,,8} is taken over all the partitions of unity $D$ in the algebra $C(X)$. If
we have $\mu(g)  = 0$ for a certain function $g\in D$, then we set the corresponding summand
in~\eqref{5,,9} to be zero independently of the value~$m(A^ng)$. And if there exists a function
$g\in D$ such that  $A^ng\equiv  0$ and simultaneously $\mu(g)>0$, then we set $\tau(\mu)  =
-\infty$.
\end{defn}

\begin{rem}\label{5..4}
As we know the Kolmogorov-Sinai definition of entropy $h(\mu)$ is based on the notion of  entropy
of a measure $\mu$ with respect to a partition of space and  it characterizes  the behavior of the
entropy of partition under its refinement by means of the mapping $ \alpha$. Recall that given a
partition of the space $X$ by measurable subsets
$$
  X =\coprod_i D_i
$$
its entropy is defined by the expression
\begin{equation}\label{5,,10}
 - \sum_i  \mu(D_i) \ln\mu(D_i).
\end{equation}
Therefore it seems reasonable to look for lower estimates for the spectral potential that contain
expressions similar to~(\ref{5,,10}).

Similarly to (\ref{5,,10}), given a partition of unity $D$ one can consider the sum
$$
 - \sum_i  \mu(g_i)\ln\mu(g_i)
$$
as the entropy of the measure $\mu$ assigned to this partition, and the sum
$$
 - \sum_i \mu(g_i)\ln m(g_i)
$$
can be considered as the corresponding \emph{relative} entropy of the measure $m$.

Then the expression
$$
 - \sum_i  \mu(g_i) \ln \frac{m(g_i)}{\mu(g_i)}
 = - \sum_i  \mu(g_i) \ln \mu(g_i) +\sum_i  \mu(g_i) \ln m(g_i)
$$
can be interpreted as the \emph{entropy distance} (difference) from $m$ to $\mu$.

Thus the expression
$$
 \inf_{m\in A^{*n}(M)} \left(-\sum_{g\in D}\mu(g) \ln\frac{m(g)}{\mu(g)}\right)
$$

\medskip\noindent
may be treated as the entropy distance from measure  $\mu$ to the image of the adjoint oper\-ator
$A^{*n}$ or the `shortage  of the entropy'. Here  the images $ L_n = A^{*n}(M)$ form a decreasing
chain, thus the distance increases and $\tau(\mu)$ characterizes the speed of the increase of the
distance.

From this point of view the case when  $ \alpha$ is reversible is fairly vivid: $ L_n = M$ and the
distance is zero.
\end{rem}

In the notation of Definition \ref{5..3} inequality  \eqref{5,,7} takes the form
$\lambda(\varphi)\ge \mu(\varphi) +\tau(\mu)$ and therefore Corollary \ref{5..2} can be rewritten
as

\begin{theorem}\label{5..5}
{\bf [The lower estimate of the spectral potential]\,}
\begin{equation}\label{5,,11}
 \lambda(\varphi)\ge \sup_{\mu\in M_\alpha}\bigl(\mu(\varphi) +\tau(\mu)\bigr).
\end{equation}
\end{theorem}

Now we are ready to formulate the principal result of the article. Its essence is that inequality
\eqref{5,,11} is in fact equality, and the restriction of $t$-entropy $\tau(\mu)$ onto the set of
invariant measures  $M_\alpha$ coincides with the dual entropy $S(\mu)$. Namely, the following
theorem holds true.

\begin{theorem}\label{5..6}
{\bf [Variational Principle in the model setting]\,} Let\/ $A\!:C(X)\to C(X)$ be a transfer
operator for a continuous mapping\/  $\alpha\!:X\to X$ of a Hausdorff compact space\/~$X$.  Then
its spectral potential\/~$\lambda(\varphi)$ satisfies the variational principle
\begin{equation}\label{5,,12}
\lambda(\varphi)= \max_{\mu\in M_\alpha} \bigl(\mu(\varphi)+\tau(\mu)\bigr),\qquad\varphi\in C(X),
\end{equation}
and\/ $t$-entropy satisfies the equality
\begin{equation}\label{5,,13}
\tau(\mu) =\inf_{\varphi\in C(X)}\bigl(\lambda(\varphi) -\mu(\varphi)\bigr),\qquad \mu\in M_\alpha.
\end{equation}
\end{theorem}

\begin{rem}\label{5..7}
 1) In fact the equalities established in  this theorem are much deeper
than simply the explicit calculation of the Legendre dual objects arising in the procedure of
spectral radius evaluation. Much more important (from our point of view) is the observation that
formula (\ref{5,,12})    links the spectral characteristics of the transfer operator  (the
left-hand part) with the stochastic characteristics ($\tau(\mu)$ in the right-hand part) of the
dynamical system.   This ideology will be developed further in  Section \ref{9..} where in
particular the interrelation between $\tau (\mu)$ and the distribution of empirical measures is
described. Namely this interrelation will serve as a key instrument in the proof of the Variational
Principle.

2) Formula (\ref{5,,12}) reveals the partition  of the process of calculation of the spectral
radius into the \emph{static component} (the first summand in the right-hand part depends only on
the weight $\varphi $)  and the  \emph{dynamical component} (the second summand depends only on the
shift $\alpha$ and the transfer operator $A$).

3) The duality established in Theorem \ref{5..6} and the thermodynamic formalism developed in
\cite{ABLS} leads naturally to introduction of the thermodynamic `ideology' into the spectral
analysis of transfer  operators. Having in mind this motivation it is reasonable to call the
functionals (measures) $\mu$ at which the maximum in the right-hand part  of~(\ref{5,,12})   is
attained the \emph{equilibrium states}. We recall in this connection that in accordance with  a
common physical point of view the equilibrium states are the states at which the system `exists in
reality'. From this point of view the duality principle adds dialectics to the spectral analysis of
transfer operators: since $\tau(\mu)$ describes the  measure of the `most typical' trajectories
(see, in particular,  Section \ref{9..}) and the value $\mu(\varphi)$ calculates the `living
con\-di\-tions' (recall the corresponding discussion in \cite{ABL-2005}) then the duality principle
tells us that the process realizes at a state having the best combination of these components.
\end{rem}

Theorem \ref{5..6} can be extended up to the case of not only positive but also nonnegative
coefficients of transfer operators. Namely, let us consider  $\psi :=e^\varphi$ and denote by
$\ell(\psi)$ the logarithm of the spectral radius of the operator $A\psi =A(\psi\,\cdot\,)$. We
have by definition that  $\ell(\psi) =\lambda(\ln\psi)$, and formula~\eqref{5,,12} takes the form
\begin{equation}\label{5,,14}
\ell(\psi)  = \max_{\mu\in M_\alpha(\cal C)} \biggl(\int_X \ln\psi\,d\mu +\tau(\mu)\biggr).
\end{equation}
Theorem \ref{5..6} tells that this formula holds true for all strictly positive functions $\psi\in
C(X)$.

\begin{theorem}\label{5..8}
Under the conditions of Theorem \ref{5..6} formula\/ \eqref{5,,14} holds true for all nonnegative
functions\/~$\psi\in C(X)$.
\end{theorem}

Clearly equality  \eqref{5,,14} is stronger than  \eqref{5,,12}. Nevertheless, it seems worthwhile
to consider formula~\eqref{5,,12} as the main form of the Variational Principle since it is more
natural from the point of view of the Legendre duality.

\section{Variational Principle: {\boldmath$C^*$}-algebraic picture}\label{6..}

We have formulated the Variational Principle for the spectral potential of a transfer operator in
the model setting when the phase space $X$ is compact,  the mapping $\alpha$ is continuous and the
transfer operator acts in  $C(X)$. In fact the most natural operator picture can be written by
means of $C^*$-algebraic language. In this section we present this picture. The $C^*$-algebraic
language, apart from its universality, will provide us with an essential extension of the range of
applications.

\begin{I}\label{6..1}
We will call $\cal C$ a \emph{base algebra} if it is a selfadjoint part of a certain commutative
\hbox{$C^*$-al}\-ge\-bra with an identity $\bf 1$. This means that there exists a commutative
$C^*$-algebra $\cal B$ with an identity $\bf 1$ such that
$$
{\cal C} = \{\,b\in {\cal B}\mid b^* = b\,\}.
$$

As is known the Gelfand transform establishes an isomorphism between  $\cal C$ and  the algebra
$C(X)$ of continuous real-valued functions on a Hausdorff compact space $X$, which is the maximal
ideal space of the algebra $\cal C$. Throughout the article we identify $\cal C$ with $C(X)$
mentioned above.
\end{I}

The next known result (see, for example, \cite{end-R}) establishes a correspondence between
endomorphisms of base algebras and dynamical systems.

\begin{theorem}\label{6..2}
If\/ $\delta\!: {\cal C} \to {\cal C}$ is an endomorphism of a base  algebra\/ $\cal C$ then there
exists an   open-closed subset\/  $Y\subset X$ and a continuous mapping\/  $\alpha\!:Y\to X$
(both\/ $Y$ and\/ $\alpha$ are uniquely defined) such that
$$
\bigl[\delta f\bigr](x) = \chi_{Y}(x)f(\alpha (x)), \qquad f\in {\cal C}, \quad x\in X,
$$
where\/ $\chi_{Y}$ is the index function of\/ $Y$. In particular if\/ $\delta({\mathbf 1}) =
{\mathbf 1}$ then\/ $Y = X$ and
\begin{equation}\label{6,,1}
\bigl[\delta f\bigr](x) = f(\alpha (x)).
\end{equation}
\end{theorem}

\begin{rem}\label{6..3}
It is clear that any endomorphism of a $C^*$-algebra $\cal B$ is completely defined by its
restriction onto the selfadjoint part $\cal C$ of $\cal B$ and  on the other hand any endomorphism
of $\cal C$ extends uniquely up to an endomorphism of $\cal B$. Therefore the  cor\-res\-pond\-ence
between endomorphisms and dynamical systems presented in the theorem can be equally described in
terms of endomorphisms of~$\cal B$.
\end{rem}

\begin{I}\label{6..4}
In what follows the pair $({\cal C},\delta)$, where $\cal C$ is a base algebra and $\delta$ is its
certain endomorphism such that $\delta({\mathbf 1}) = {\mathbf 1}$, will be called a
$C^*$-\emph{dynamical system}, and the pair  $(X, \alpha)$ described in Theorem \ref{6..2} will be
called the \emph{dynamical system corresponding to} $({\cal C},\delta)$. The algebra ${\cal C}$
will be also called the \emph{base algebra of the dynamical system} $(X,\alpha)$.

Throughout the paper  notation  ${\cal C}$, $\delta$, $X$, $\alpha$ will denote the objects
introduced above and we will use either of them (say $\delta$ or $\alpha$) for convenience reasons
(once $\alpha $ is chosen then $\delta$ is defined uniquely by (\ref{6,,1}) and vice versa).
\end{I}

Definition \ref{2..1} of transfer operator can be  naturally rewritten in terms of $C^*$-dynamical
systems.

\begin{defn}\label{6..5}
Let $({\cal C},\delta)$ be a $C^*$-dynamical system. A linear operator $A\!:\cal C\to \cal C$ will
be called a \emph{transfer operator}, if it possesses the following two properties

a) $A$ is positive (it maps nonnegative elements of $\cal C$ into nonnegative ones);

b) it satisfies the  \emph{homological identity}
\begin{equation}\label{6,,3}
 A\bigl((\delta f)g\bigr) =fAg \quad \textrm{for all}\ \ f,g\in \cal C.
\end{equation}
If in addition this operator maps $\mathbf 1$ into $\mathbf 1$ we will call it a \emph{conditional
expectation operator.}
\end{defn}

\begin{rem}\label{6..6}
Any transfer operator $A\!:\cal C\to \cal C$ can be naturally extended up to a transfer operator on
${\cal B} = {\cal C} + i{\cal C}$ by means of the formula
$$
A (f +ig)= Af + iAg.
$$
On the other hand given any transfer operator on $\cal B$, its restriction to $\cal C$ (which is
well defined in view of property a) of Definition \ref{6..5}) is also a transfer operator.
Therefore transfer operators can be equivalently introduced as by means of  $C^*$-algebra $\cal B$
so also by means of its selfadjoint part --- the base algebra $\cal C$. We prefer to exploit the
base algebra since in what follows we use the Legendre transform which is an essentially
real-valued object.
\end{rem}

\begin{rem}\label{6..7}
 1) In view of (\ref{6,,3}) the conditional expectation operator satisfies
the equality
$$
 A(\delta f) =A\bigl((\delta f){\mathbf 1}\bigr) =fA{\mathbf 1}= f\quad
 \textrm{for every}\ \  f\in \cal C.
$$
Thus this operator is a positive left inverse to the mapping $f \mapsto \delta f$.

2) If $A\!: {\cal C} \to {\cal C}$ is a transfer operator and  $\frac{1}{A\mathbf 1}\in\cal C$ then
$\frac{1}{A\mathbf 1}A$ is a conditional expectation operator.
\end{rem}

A more detailed analysis and various types of $C^*$-dynamical systems and transfer operators will
be presented  in the next Section \ref{7..}.

\begin{I}\label{6..8}
Let $(\cal C,\delta)$ be a $\cal C^*$-dynamical system and $(X,\alpha)$ be the corresponding
dynamical system. We denote by  $M(\cal C)$ the set of all positive normalized linear functionals
on~$\cal C$ (which take nonnegative values on nonnegative elements and  are equal to $1$ on the
unit). Since we are identifying $\cal C$ and $ C(X)$,  the Riesz theorem implies that the set
$M(\cal C)$ can be identified with the set of all regular Borel probability measures on~$X$ and the
identification is established by means of the formula
\begin{equation*}
 \mu(\varphi) = \int_X \varphi\, d\mu,  \qquad \varphi \in{\cal C}=C(X),
\end{equation*}

\smallskip\noindent
where $\mu$ in the right-hand part is a measure on $X$ assigned to the functional $\mu\in M (\cal
C)$ in the left-hand part. That is why with a slight abuse of language we will call elements of
$M(\cal C)$ measures.

A measure  $\mu\in M(\cal C)$ is called  \emph{$\delta$-in\-va\-ri\-ant} if for each $f\in \cal C$
we have $\mu(f)=\mu(\delta f)$. The set of all $\delta$-in\-va\-ri\-ant measures from  $M(\cal C)$
will be denoted by $M_\delta(\cal C)$. Clearly, in terms of the dynamical system $(X,\alpha)$ the
condition $\mu(f) =\mu(\delta f)$ is equivalent to the condition $\mu(f) =\mu(f\circ\alpha)$,\ \
$f\in C(X)$. Therefore $M_\delta(\cal C)$ can be identified with the set of all
$\alpha$-in\-vari\-ant Borel probability measures on X.
\end{I}

\begin{I}\label{6..9}
Let $A\!: \cal C\to \cal C$ be a fixed transfer operator for a $C^*$-dynamical system~$({\cal C},
\delta)$. Just as in Section \ref{2..} we define the family of operators $A_\varphi\!: \cal C\to
\cal C$, where $\varphi \in {\cal C}$, by means of the formula $A_\varphi f := A(e^\varphi f)$.
Evidently, all the operators of this family are transfer operators for $({\cal C}, \delta)$ as
well.
\end{I}

In view of the foregoing observations the definitions of the spectral potential $\lambda(\varphi)$
(equality (\ref{2,,2})) and  $t$-entropy (Definition \ref{5..3}) can be  carried over  word by word
to the $C^*$-dynamical systems case. Here they are:
\begin{gather}\label{6,,5}
 \lambda(\varphi)  = \lim_{n\to\infty}\frac{1}{n}\ln \norm{A_\varphi^n\mathbf 1},\qquad \varphi\in
 \cal C;\\[9pt]\label{6,,6}
 \tau(\mu)  := \inf_{n\in\mathbb N}\frac{\tau_n(\mu)}{n}\,,\qquad
 \tau_n(\mu)  := \inf_D\tau_n(\mu,D),\\[6pt]\label{6,,7}
 \tau_n(\mu,D)  := \sup_{m\in M(\cal C)}\sum_{g\in D}\mu(g) \ln\frac{m(A^ng)}{\mu(g)}\,,\qquad
 \mu\in M_\delta(\cal C).
\end{gather}
The infimum in \eqref{6,,7} is taken over all the partitions of unity $D$ in the algebra $\cal C$.

Once this is done Theorem  \ref{5..6} takes the following form.

\begin{theorem}\label{6..10}
{\bf [Variational Principle: \boldmath{$C^*$}-formulation]\,} Let\/ $({\cal C},\delta)$ be a\/
$C^*$-dyna\-mic\-al system,  $A\!:\cal C\to \cal C$ be a certain trans\-fer operator for\/ $({\cal
C},\delta)$, and\/ $A_\varphi =A(e^\varphi\,\cdot\,)$ for all\/ $\varphi\in\cal C$. Then the
corresponding spectral potential\/ $\lambda(\varphi)$ satisfies the variational principle
\begin{equation}\label{6,,4}
 \lambda(\varphi)  = \max_{\mu\in M_\delta(\cal C)} \bigl(\mu(\varphi) +\tau(\mu)\bigr), \qquad
 \varphi\in \cal C,
\end{equation}

\smallskip\noindent where $M_\delta(\cal C)$ is the set of all positive normalized
$\delta$-in\-va\-ri\-ant linear functionals on\/ $\cal C$.
\end{theorem}

The proof of Theorem \ref{6..10}   exploits  a number of special properties of $t$-entropy.
Therefore  before proceeding to the proof of the Variational Principle itself (which will be given
in Section \ref{10..}) we will implement the necessary analysis of $t$-entropy. The results of this
analysis are presented in two Sections \ref{8..} and \ref{9..} where in the first one we describe
the analytic properties of $t$-entropy while the second one is devoted to its statistical
properties.

The $C^*$-analogue to Theorem \ref{5..8} will be formulated and proved in Section \ref{11..}.

\section[Types  of \boldmath{$C^*$}-dynamical systems and transfer operators]{Types
 of \boldmath{$C^*$}-dynamical systems\\ and transfer operators}\label{7..}

In this section we present a number of types (which we call examples) of $C^*$-dynamical systems
and transfer operators that show, in particular, how far away from the model situation presented in
Section \ref{2..} one can move.   In addition we give a complete description of transfer operators.

Let us start with  certain examples of $C^*$-dynamical systems and base algebras.

\begin{exm}
\label{7..1} \em Let $Y$ be a measurable space with a  $\sigma$-algebra $\goth A$ and $\beta\!:Y\to
Y$ be a measurable mapping. We denote by  $(Y,\beta)$ the discrete time dynamical system
generated by the mapping $\beta$ on the phase space $Y$. Let $\cal B$ be any Banach algebra such
that

a) $\cal B$ consists of bounded real-valued measurable functions on  $Y$,

b) it is supplied with the uniform norm,

c) it contains the unit function, and

d) it is $\beta$-invariant (that is $f\circ\beta\in\cal B$ for all  $f\in\cal B$).

Clearly the mapping $\delta\!: \cal B\to \cal B$ given by $\delta(f) := f\circ\beta$ is an
endomorphism of $\cal B$ and therefore $({\cal B}, \delta)$ is a $C^*$-dynamical system with the
base algebra $\cal B$.
\end{exm}

\begin{exm}\label{7..2}
As a particular case of the base algebra  in the previous example one can take the
algebra  of all bounded real-valued measurable functions on  $Y$. We will denote this algebra
by~$B(Y)$.
\end{exm}

\begin{exm}\label{7..3}
Let $(Y,\goth A,m)$ be a measurable space with a probability measure $m$, and let~$\beta$ be a
measurable mapping such that $m\bigl(\beta^{-1}(G)\bigr)\le Cm(G),\ \ G\in \goth A$, where the
constant $C$ does not depend on $G$. In this case one can take as a base algebra the space
$L^\infty(Y,m)$ of all essentially bounded real-valued measurable functions on $Y$ with the
essential supremum norm.
\end{exm}

\begin{rem}\label{7..4}
1) If, as  in Example \ref{7..2},\ \ $\cal C =B(Y)$ then the elements of~$M(\cal C)$ can be
naturally identified with finitely-additive probability measures on the $\sigma$-algebra~$\goth A$
by means of the equality $\mu(f) =\int_Y f\,d\mu,\ \ f\in B(Y)$.

2) If, as in Example \ref{7..3},\ \ $\cal C = L^\infty (Y,m)$ then $M(\cal C)$ consists of
finitely-additive probability measures on $\goth A$ which are absolutely continuous with respect
to~$m$ (that is they are equal to zero on the sets of zero measure $m$).

3) In  Example \ref{7..3} the set $M_\alpha(\cal C)$ is the subset of $M(\cal C)$ consisting of
measures $\mu$ such that $\mu(\beta^{-1}(G)) =\mu(G)$ for each measurable set~$G$.

4) It should be emphasized  that in general given a concrete functional algebra its endomorphism is
\emph{not} necessarily generated by \emph{a point mapping} of the domain. For example, if $\cal C =
L^\infty (Y,m)$ then its endomorphisms are generated by \emph{set mappings}, that do not `feel'
sets of measure zero
(see, for example \cite{Walt2},~Chapter~2). Thus not every endomorphism of \ $
L^\infty (Y,m)$ is generated by a certain measurable mapping~$\beta$ as in Example  \ref{7..3}.

On the other hand, as Theorem \ref{6..2} tells,  on the maximal ideals level any endomorphism is
induced by a certain point mapping. Therefore raising the apparatus of investigation  to the
$C^*$-algebraic level (Definition \ref{6..5}) we not only essentially extend the sphere of the
objects under consideration but additionally  can always exploit   point mappings in the  study of
transfer operators independently of their concrete origin (see in this connection the general
description of transfer operators given  below).
\end{rem}

Now let us present  certain  examples of transfer operators.

\begin{exm}\label{7..5}
Let $X$ be a compact space, $\alpha\!: X\to X$ be a homeomorphism and ${\cal C}= C(X)$. Consider a
weighted shift operator $A :C(X)\to C(X)$ given by the formula
$$
\bigl[Af\bigr](x) = a(x)f(\alpha^{-1}(x)),
$$
where $a\in C(X)$ is a certain nonnegative function. Evidently $A$ is a transfer operator.
\end{exm}

The next example can be considered as the main model example for transfer operators discussed in
the article.

\begin{exm}\label{7..6}
Let $(Y,\goth A)$ be a measurable space with a  $\sigma$-finite measure  $m$, and let $\beta $ be a
measurable mapping such that for all measurable sets  $G\in\goth A $ the following estimate holds
$$
m\bigl(\beta^{-1}(G)\bigr) \le Cm(G),
$$
where the constant $C$ does not depend on $G$. For example, if the measure~$m$ is $\beta$-invariant
one can set~$C=1$. Let us consider the space $L^1(Y,m)$ of real-valued integrable functions and the
shift operator that takes every function $f\in L^1(Y,m)$ to~$f\circ\beta$. Clearly the norm of this
operator does not exceed~$C$. The mapping $\delta f :=f\circ\beta$ acts also on the space
$L^\infty(Y,m)$ and it is  an endomorphism of this space. As is known, the dual space to $L^1(Y,m)$
coincides with $L^\infty(Y,m)$. Define the linear operator $A\!: L^\infty(Y,m) \to L^\infty(Y,m)$
by the identity
\begin{equation*}
\int_Y f\cdot g\circ\beta\,dm \,\equiv\, \int_Y (Af)g\,dm, \qquad g\in L^1(Y,m).
\end{equation*}
In other words  $A$ is the adjoint operator to the shift operator in $L^1(Y,m)$. If one takes as
$g$ the index functions of measurable sets  $G\subset Y$, then the latter identity takes the form
$$
\int_{\beta^{-1}(G)}f\,dm \equiv \int_G Af\,dm .
$$
Therefore  $Af$ is nothing else than the Radon--Nikodim density of the additive set function
$\mu_f(G) =\int_{\beta^{-1}(G)}f\,dm$. Evidently, the operator  $A$ is positive and satisfies the
homological identity
$$
A\bigl((\delta f)g\bigr) =fAg, \qquad f,g\in L^\infty(X,m).
$$

\medskip\noindent
We see that $A$ is a  transfer operator (for the $C^*$-dynamical system $(L^\infty(Y,m), \delta)$).
And in the case when $m$ is   $\beta$-invariant measure   it is a conditional expectation operator.
\end{exm}

\noindent{\bf Transfer operators and positive functionals.} By identifying the base algebra $\cal
C$ with $C(X)$ one can also obtain  a certain `more explicit' description of transfer operators
linking them with special families of positive functionals.  Here it is.

Let, as usual, $X$ be a compact space, $\alpha\!: X\to X$ be a continuous mapping, and $A\!:
C(X)\to C(X)$ be a certain transfer operator.

For every point  $x\in X$ define the functional $\phi_x$ according to the formula
\begin{equation}\label{7,,1}
 \phi_x (f) := \bigl[{A}f\bigr](x), \qquad f\in C(X).
\end{equation}
Evidently, $\phi_x$ is a positive functional.

There are two possibilities for  $x$.

1) $[A{\mathbf 1}] (x)=0$. This means that $\phi_x ({\mathbf 1})=0$ which implies $\phi_x =0$ due
to the positivity of $\phi_x$.

2) $[A{\mathbf 1}](x)\neq 0$. In this case $\phi_x \neq 0$ and $\phi_x $ defines a certain measure
$ \nu_x$  on $X$.

The homological identity implies also that for any $f \in C(X)$ we have
$$
 \bigl[A(f\circ\alpha)\bigr](x) =\bigl[A(f\circ\alpha\cdot{\mathbf 1})\bigr](x) = f(x)\cdot
 A{\mathbf 1}(x),
$$
and therefore
$$
 \frac{1}{A{\mathbf 1}(x)}\,\phi_x (f \circ \alpha)= f (x),
$$
which means that
\begin{equation}\label{7,,2}
 \mathop{\rm supp} \nu_x \subset \alpha^{-1} (x).
\end{equation}

\medskip

Clearly, the mapping  $x \to \phi_x$ is $^*$-weakly continuous on  $X$.

Observe also that if $x \notin \alpha (X)$ then $A{\mathbf 1}(x)=0$. Indeed, if $A{\mathbf
1}(x)\neq 0$ then choosing a function $f\in C(X)$ such that
$$
 f|_{\alpha (X)} =0\quad \textrm{and} \quad f(x)=1
$$
and exploiting   the homological identity one obtains
$$
 0= \frac{1}{{A\mathbf 1}(x)}\bigl[A(f\circ\alpha)\bigr](x) = f(x)=1
$$
thus arriving at a contradiction.

The  objects  presented above in fact give a complete description of transfer operators in $C(X)$
since one can easily verify  that every $^*$-weakly  continuous mapping  $x \mapsto \phi_x$, where
$\phi_x$ are positive functionals such that

a) \  $\phi_x =0,\ \ x\notin \alpha (X)$,

b) \ $\phi_x$ satisfies (\ref{7,,2}) for $x\in \alpha (X)$ (here $\phi_x$ may be $0$ as well)

\noindent
defines a certain transfer operator $A\!:C(X)\to C(X)$ acting according to formula~(\ref{7,,1}).

\begin{rem} \label{7..7}
1) If $\alpha \!: X \to X$ is a homeomorphism then the foregoing discussion shows that any transfer
operator $A\!: C(X) \to C(X)$ is a weighted shift operator mentioned in Example \ref{7..5}.

2) In general for a continuous mapping $\alpha\!:X\to X$ even a conditional expectation in~$C(X)$
(if it exists) is not defined in a unique way.

For example, consider $X = \mathbb T^1 = \mathbb R/\mathbb Z$ and the mapping $\alpha\!:X\to X$
defined by the formula $\alpha(x) = 2x\pmod{1}$. Take any continuous function $\rho$ on $X$ having
the properties
$$
0\le \rho (x) \le 1,\quad \rho\left(x+\frac{1}{2}\right) + \rho (x)=1,\qquad x\in X,
$$
Evidently, the operator
$$
 \bigl[Af\bigr](x) :=  \sum_{y\in\alpha^{-1}(x)}f(y)\rho(y) = f\left(\frac{x}{2}\right)
\rho\left(\frac{x}{2}\right) + f\left(\frac{x+1}{2}\right)\rho\left(\frac{x+1}{2}\right)
$$
is a conditional expectation operator for $(X,\alpha)$.

3) Recalling Remark \ref{7..4} \, 4) \, we have to stress that in general given a concrete
functional algebra and its  endomorphism then a transfer operator is \emph{not} necessarily
associated with \emph{a point mapping} of the domain.
\end{rem}

\section{Properties of $\boldsymbol{t}$-entropy }\label{8..}

In this section we prove a number of properties of  $t$-entropy and in particular its upper
semicontinuity.

Let us start with a simple but important observation.

\begin{prop}\label{8..1}
The functional\/ $\tau(\mu)$ satisfies the inequality\/ $\tau(\mu)\le \lambda(0)$.
\end{prop}

\proof. Since logarithm is a concave function it follows that for any partition of unity $D$ in
\eqref{6,,7} we have
\begin{align*}
 \tau_n(\mu,D) &=\sup_{m\in M(\cal C)}\sum_{g\in D}\mu(g) \ln\frac{m(A^ng)}{\mu(g)} \le \sup_{m\in
 M(\cal C)}\ln\sum_{g\in D} m(A^ng)\\[3pt]
 &= \sup_{m\in M(\cal C)}\ln m(A^n\mathbf 1) = \ln\norm{A^n}.
\end{align*}
Which implies the desired inequality.\qed

\begin{rem}\label{8..2}
In particular, if  $A$ is a conditional expectation operator (i.\,e.\ $A\mathbf 1 =\mathbf 1$) then
$\norm{A^n} =1$ and therefore  $\lambda (0)=0$. Thus in this case $\tau(\mu)\le 0$.
\end{rem}

\begin{prop}\label{8..3}
If $A\!:{\cal C} \to {\cal C}$ is an invertible conditional expectation operator then\/
$\tau(\mu)=0$ for  any\/ $\mu\in M({\cal C})$.
\end{prop}

\proof. Since $A$ is a conditional expectation operator we have $ A\bigl(\delta f) = f,\ \ f\in
\cal C$, and hence $A^{-1}f = \delta f$. This implies in particular that $A^{-1}\!: \cal C\to \cal
C$ is a positive mapping. For any $\mu\in M({\cal C})$ and $n\in \mathbb N$ define the functional
$m$ by the formula
$$
m(f) := \mu (A^{-n} (f)), \qquad f\in \cal C.
$$
Clearly, $m\in M ({\cal C})$. For this $m$ we have
$$
 \sum_{g\in D}\mu(g) \ln\frac{m(A^ng)}{\mu(g)} =0,
$$
and therefore $\tau_n(\mu,D)\ge 0$ and $\tau(\mu) \ge 0$. Combining the latter inequality with
Remark~\ref{8..2} we obtain $\tau(\mu) = 0$. \qed

\begin{prop}\label{8..4}
The functionals\/ $\tau_n(\mu,D)$ and\/ $\tau(\mu)$ depend concavely on\/ $\mu\in M(\cal C)$.
\end{prop}

\proof. Let $\mu_1,\mu_2\in M(\cal C)$ and $\mu  = p_1\mu_1+p_2\mu_2$, where  $p_1+p_2  = 1$ and
$p_1,p_2\ge 0$. Then for any element $g$ of a partition of unity~$D$ we have
\begin{equation*}
p_1\mu_1(g)\ln\frac{m_1(A^n g)}{\mu_1(g)} + p_2\mu_2(g)\ln\frac{m_2(A^n g)}{\mu_2(g)} \le
\mu(g)\ln\frac{[p_1m_1+p_2m_2](A^n g)}{\mu(g)}.
\end{equation*}
Let us sum this inequality over $g\in D$ and pass to the supremums with respect to  $m_1$ and
$m_2$. As a result we obtain the inequality
$$
p_1\tau_n(\mu_1,D) +p_2\tau_n(\mu_2,D) \le \tau_n(\mu,D).
$$
It means that the functional  $\tau_n (\mu,D)$ is concave with respect to $\mu$. This along with
\eqref{5,,8} implies the concavity of  $\tau(\mu)$. \qed \proofskip

\medskip

Let  $D =\{g_1,\dots,g_k\}$ be a partition of unity in a base algebra~$\cal C$. We eliminate from
this partition the elements~$g_i$ such that $A^ng_i = 0$ and put $D'_n := \{\,g\in D\mid A^ng \ne
0\,\}$. Denote by $M(D)$ and $M(D'_n)$ the sets of all probability measures on finite sets $D$
and~$D'_n$, respectively. Clearly, these $M(D)$ and $M(D'_n)$ are finite-dimensional simplexes. If
one extends each measure $\mu\in M(D'_n)$ by zero to $D\setminus D'_n$, then the simplex  $M(D'_n)$
becomes a certain face of  $M(D)$. Note that formula \eqref{6,,7} defines the functions
$\tau_n(\mu,D)$ not only for measures  $\mu\in M(\cal C)$ but also for measures  $\mu\in M(D)$.

\begin{prop}\label{8..5}
The function\/ $\tau_n(\,\cdot\,,D)$ is continuous on the set\/ $M(D'_n)$ and it is equal
to\/~$-\infty$ on\/~$M(D)\setminus M(D'_n)$.
\end{prop}

\proof. Let us take  a measure  $m'\in M(\cal C)$ such that  $m'(A^ng) >0$ for all  $g\in D'_n$.
Fix a (small) positive number  $\varepsilon$. For every measure  $m\in M(\cal C)$ consider the new
measure
\begin{equation}\label{8,,1}
 m_\eps := (1-\eps)m +\eps m'.
\end{equation}
Then for all $g\in D'_n$ the following inequalities hold
$$
m_\eps(A^ng)\ge \eps m'(A^ng)>0
$$
and, on the other hand,
$$
 m_\eps(A^ng)\le \norm{A^n\mathbf 1}.
$$
These inequalities imply the existence of a (large) number $C(\eps)$ such that
$$
 \bigl|\ln m_\eps(A^ng)\bigr|\le C(\eps)\quad \text{for all}\ \ m\in M(\cal C),\,\ g\in
D'_n.
$$
Let us introduce the notation
\begin{gather}\label{8,,2}
 \psi(\mu,m)  := \sum_{g\in D'_n} \mu(g)\ln m(A^ng),\\[6pt]\label{8,,3}
 \psi_\eps(\mu) := \sup_{m\in M(\cal C)}\psi(\mu,m_\eps).
\end{gather}
Then for all measures  $\mu,\nu\in M(D'_n)$ we have
\begin{equation*}
\bigl|\psi(\mu,m_\eps) -\psi(\nu,m_\eps)\bigr|\le C(\eps) \sum_{g\in D'_n}\bigl|\mu(g)
-\nu(g)\bigr|.
\end{equation*}
Therefore the function $\psi_\eps(\mu)$ depends continuously on $\mu\in M(D'_n)$.

Note that  (\ref{8,,1}) and (\ref{8,,3}) imply
\begin{equation}\label{8,,4}
 \psi_0(\mu)\ge \psi_\eps(\mu).
\end{equation}
On the other hand since $m_\eps \ge (1-\eps)m$ we have
\begin{equation} \label{8,,5}
  \psi_\eps(\mu)\ge \psi_0(\mu)+ \ln(1-\eps).
\end{equation}
In view of (\ref{8,,4}) and (\ref{8,,5}) the function $\psi_0(\mu)$ is the  uniform limit of the
functions $\psi_\eps(\mu)$ as $\eps\to 0$. Thus it is continuous as well. Finally, the difference
of two continuous functions
$$
 \psi_0(\mu)\,-\,\sum_{g\in D'_n} \mu(g)\ln \mu(g)
$$
coincides with  $\tau_n(\mu,D)$ and so it is continuous on the set~$M(D'_n)$.

The second part of the proposition follows from~\eqref{5,,9}. \qed

\begin{prop}\label{8..6}
The functional\/ $\tau(\mu)$ is upper semicontinuous on\/ $M(\cal C)$.
\end{prop}

\proof. By Proposition \ref{8..5} the function  $\tau_n(\mu,D)$ is upper semicontinuous (with
respect to $^*$-weak topology) on $M(\cal C)$. Therefore, $t$-entropy
$$
 \tau(\mu) = \inf_{n,D}\tau_n(\mu,D)/n.
$$
also possesses this property. \qed \proofskip

The definition of $\tau(\mu)$ is rather complicated and its reduction to a simpler form in general
is problematic. Let us present an example when the corresponding calculation leads us to simpler
expressions.

It is also worth remarking that the explicit calculation of  $ \tau(\mu)$ for a concrete measure
gives us a lower estimate for the spectral potential.

\begin{exm}\label{8..7}
Throughout this example we confine ourselves to the case of Perron--Fro\-be\-nius
operator~\eqref{1,,3} when the preimage of every point consists of finite number of points.

First, let $ x_0 $ be a fixed point of the mapping $\alpha$. Then the measure $\mu =\delta_{x_0}$
is a probability invariant measure.

Let $D = \{g_0, g_1, \dots , g_k\}$ be a partition of unity. Then $\mu(g_i) = g_i(x_0)$. Since
in~(\ref{5,,8}) we take the infimum over all partitions of unity it follows that one has to
consider the partitions for which $g_0(x_0) =1 $ and $g_i(x_0) =0 $ for $ i\ne 0$. Then the sum in
expression  (\ref{5,,9}) for  $\tau_n(\mu,D)$ reduces to a single summand
\begin{equation}\label{8,,6}
 \tau_n(\mu,D) = \sup_{m\in M} \ln{m(A^ng_0)}.
\end{equation}
For $m = \mu =\delta_{x_0}$ we have $m(A^ng_0) =A^ng_0(x_0)$. By assumption the operator  $A$ acts
according to  formula
$$
 Ag(x) =\sum_{y\in\alpha^{-1}(x)}\psi(y)g(y),
$$
where $ \psi \in C(X)$ is a fixed nonnegative function, and thus
$$
 A^ng(x_0) =\sum_{y\in\alpha^{-n}(x_0)}\left[\prod_{j =0}^{n-1}\psi(\alpha^j(y))\right]g(y).
$$
Given $n$ and $x_0$, one can choose a partition of unity such that
$g_0(y) = 0$ for $ y\in\alpha^{-n}(x_0)$,\ \ $y\ne x_0$. Then $A^ng_0(x_0) =\psi (x_0)^n$ and therefore
$$
 \tau_n(\mu,D)\geq  n\ln\psi(x_0)
$$
and
$$
\tau(\mu) \geq  \ln \psi (x_0).
$$
One can find an upper estimate for the supremum in (\ref{8,,6}) by using the following reasoning.
Given  $n$ and $\varepsilon>0$ one can choose a function  $ g_0 $ satisfying the conditions
mentioned above and having the support in a sufficiently small neighborhood of the point $x_0$.
Then $A^ng_0(x) \leq A^ng_0(x_0) +\varepsilon$ and it follows that for any probability measure $m$
one has  $ m(A^ng_0)\leq A^ng_0(x_0) +\varepsilon$. Thus in the situation considered we have
$$
 \tau(\delta_{x_0}) = \ln \psi (x_0).
$$

Now suppose $x_0$ is a periodic point with the period $N$, that is $\alpha^N(x_0) =x_0$. Then the
measure
$$
 \delta_{x_0,N} =\frac{1}{N} \sum_{j=0}^{N-1}\delta_{\alpha^j(x_0)}
$$
is an invariant probability measure. Applying the foregoing reasoning to the operator $A^N$ we
obtain
\begin{equation}\label{8,,7}
 \tau(\delta_{x_0,N}) = \frac{1}{N}\ln \prod_{j =0}^{N-1}\psi(\alpha^j(x_0)).
\end{equation}

In the wavelet theory there arises the so-called \emph{subdivision operator}, which is defined in
the following way. Let $\mathbb{T}^m = \mathbb{R}^m\!/\mathbb Z^m$ be the $m$-dimensional torus. We
consider a matrix~$M$ with integer elements such that the absolute values of all its eigenvalues
are greater than~$1$. It generates the mapping $\alpha\!:\mathbb{T}^m\to \mathbb{T}^m$ according to
the formula $\alpha(x) = Mx \pmod{\mathbb Z^m}$. This mapping preserves the Lebesgue measure and it
is an expanding map, and the number of preimages of every point is $|\mathop{\mathrm{det}}M|$.

The subdivision operator acts on the space  $ L^2(\mathbb{T}^m)$ according to the formula
$$
 S_a^M u(x) = a(x)u(\alpha(x)),
$$
where $a$ is a fixed function. It is a weighted shift operator generated by irreversible mapping.
This operator is exploited in the iteration procedure of the construction of the wavelets base. The
information on its spectral radius is required for the  convergence guarantee of the procedure
mentioned. As we have already observed  in  Introduction  the spectral radii of weighted shift and
transfer operators are tightly related (cf.~(\ref{1,,9})). Their interrelation will be discussed
further in Section \ref{13..}.

In the paper by Didenko \cite{D} a series of lower estimates for the spectral radius of the
subdivision operator was obtained. In fact, some of these estimates has been derived with the help
of the usage of periodic points of the mapping $\alpha$.

We dwell on this in such detail to emphasize the fact that such estimates follow directly from
(\ref{8,,7}). To be true, we have to note that in \cite{D} there was considered a more complicated
situation as well, namely the subdivision operators  with  \emph{matrix} coefficients in the space
of vector functions. The whole of the theory presented in this article is not applicable readily to
these operators.
\end{exm}

%
%
%
%

\section{Entropy Statistic Theorem}\label{9..}

In this section we prove a certain theorem on the statistic of distribution of empirical measures.
This result is important in its own right and plays for $\tau(\mu)$ the role similar to that the
Shannon--McMillan--Breiman theorem plays for $h(\mu)$ (see, for example \cite{Bill}, Chapter~4).
The Entropy Statistic Theorem (Theorem \ref{9..1}) not only uncovers the statistical nature of
$\tau(\mu)$ but also serves as the main technical instrument in the proof of the Variational
Principle.
\medskip

Consider a dynamical system $(X,\alpha)$, where $\alpha\!:X\to X$ is a continuous mapping of a
Hausdorff compact space $X$. Recall that by $M$ we denote the set of all Borel probability measures
on $X$. Let $x$ be an arbitrary point of $X$. The \emph{empirical measures}  $\delta_{x,n}\in M$
are defined by the formula
\begin{equation}\label{9,,1}
 \delta_{x,n}(f)  = \frac{1}{n}\Bigl(f(x) +f(\alpha(x))+\,\dotsm\,
 +f\bigl(\alpha^{n-1}(x)\bigr)\Bigr), \qquad f\in  C(X).
\end{equation}
Evidently, the measure $\delta_{x,n}$ is concentrated on the trajectory of the point $x$ of length
$n$.

We endow the set $M$ with the $\,^*$-weak topology of the dual space to $C(X)$. Given a measure
$\mu\in M$ and its certain neighborhood $O(\mu)$ we define the sequence of sets $X_n(O(\mu))$ as
follows:
\begin{equation}\label{9,,2}
 X_n(O(\mu))  := \{\,x\in X\mid \delta_{x,n}\in O(\mu)\,\}.
\end{equation}

\begin{theorem}\label{9..1}
{\bf [Entropy Statistic Theorem]\,} Let\/ $(X,\alpha)$ be a dynamical system and\/ $A\!: C(X) \to
C(X)$ be a certain transfer operator for\/ $(X,\alpha)$. Then for any measure\/ $\mu\in M$ and any
number\/ $t>\tau(\mu)$ there exist a neighborhood\/ $O(\mu)$ in the\/ $^*$-weak topology, a (large
enough) number\/~$C(t,\mu)$ and a sequence of functions\/ $\chi_n\in C(X)$ majorizing the index
functions of the sets\/~$X_n(O(\mu))$ such that for all\/ $n$ the following estimate holds
$$
 \norm{A^n\chi_n} \le C(t,\mu)e^{nt}.
$$
\end{theorem}

To prove this theorem we need a number of auxiliary results and we start with their consideration.

\medskip

Let us fix a natural number  $n$ and a partition of unity  $D$ in $C(X)$. Let, as above, $D'_n =
\{\,g\in D \mid A^ng\not\equiv 0\,\}$, and the symbols $M(D)$ and $M(D'_n)$ denote the finite
dimensional simplexes consisting of all probability measures on $D$ and $D'_n$, respectively. As it
was observed the simplex $M(D'_n)$ is naturally embedded in $M(D)$: the measures from $M(D'_n)$ can
be extended onto  $D\setminus D'_n$ by zero. Given $\mu\in M(D'_n)$ there exists a measure
$m_\mu\in M$  at which the supremum in~\eqref{5,,9} is attained. In other words,
\begin{equation}\label{9,,3}
\tau_n(\mu,D)  = \sup_{m\in M}\sum_{g\in D'_n} \mu(g)\ln\frac{m(A^ng)}{\mu(g)}  = \sum_{g\in
D'_n}\mu(g)\ln\frac{m_\mu(A^ng)}{\mu(g)}.
\end{equation}
In general the correspondence  $\mu\mapsto m_\mu$ may  be not single-valued. However, for
convenience of presentation we will assign to every measure $\mu\in M(D'_n)$ a certain
\emph{single} measure $m_\mu\in M$ satisfying the equality~\eqref{9,,3}, and thus we will assume
that there is fixed a \emph{single-valued} mapping $\mu\mapsto m_\mu$.

\begin{lemma}\label{9..2}
Given a single-valued mapping\/ $\mu\mapsto m_\mu$, satisfying equality\/ \eqref{9,,3}, every
summand\/ $\mu(g)\ln\bigl(m_\mu(A^ng)\big/\mu(g)\bigr)$ in the right-hand part of\/ \eqref{9,,3} is
a  bounded function with respect to\/ $\mu\in M(D'_n)$ and tends to zero when\/ $\mu(g) \to 0$.
\end{lemma}

\proof. First, note that  $m_\mu(A^ng)\le\norm{A^n\mathbf 1}$. Hence the expression
$\mu(g)\ln\bigl(m_\mu(A^ng)\big/ \mu(g)\bigr)$ is bounded from above, and its
$\mathop{\mathrm{sup\,lim}}$ as $\mu(g)\to 0$ is nonpositive. By Proposition \ref{8..5} the
function $\tau_n(\mu,D)$ depends continuously on $\mu\in M(D'_n)$. Thus, it is bounded on
$M(D'_n)$. So, as all the summands in the right-hand part of~\eqref{9,,3} are bounded from above,
the forgoing reasoning shows that they are bounded from below as well.

Now we finish the remaining part of the proof arguing by contradiction. Let $g_0\in D'_n$. Suppose
that the value of the expression
$$
\mu(g_0)\ln\frac{m_\mu(A^ng_0)}{\mu(g_0)\bigr)}
$$
does not tend to zero as $\mu(g_0)\to 0$. Then there exists a sequence of measures $\mu_i\in
M(D'_n)$ and a number $\eps>0$ such that
\begin{equation}\label{9,,4}
\mu_i(g_0)\to 0
\end{equation}
and at the same time
\begin{equation}\label{9,,5}
 \mu_i(g_0) \ln\frac{m_{\mu_i}(A^ng_0)}{\mu_i(g_0)\bigr)} <-\eps .
\end{equation}

\medskip\noindent
Without loss of generality one can assume that the sequence $\mu_i$ tends to a certain measure
$\nu\in M(D'_n)$. The corresponding sequence of measures $m_{\mu_i} \in M$ has at least one limit
point  $m\in M$. Passing to subsequences one can gain for every  $g\in D'_n$ the equality
$$
 \lim_{i\to\infty} m_{\mu_i}(A^ng)  = m(A^ng).
$$
Let us consider any function  $g\in D'_n$. Obviously, if $\nu(g) >0$, then
$$
 \lim_{i\to\infty}\mu_i(g) \ln m_{\mu_i} (A^ng)  = \nu(g) \ln m(A^ng);
$$
if $\nu(g) = 0$, then, as it has been already observed,
$$
\mathop{\mathrm{lim\,sup}}_{i\to\infty} \mu_i(g)\ln m_{\mu_i} (A^ng) \le 0;
$$
and if $g = g_0$, then (\ref{9,,5}) and $\mu_i(g_0)\to 0$ imply
$$
\mathop{\mathrm{lim\,sup}}_{i\to\infty} \mu_i(g_0)\ln m_{\mu_i} (A^ng_0) \le -\eps .
$$
Observe also that since $\nu(g_0)=0$ it follows that the corresponding summand in~\eqref{5,,9} for
$\tau_n(\nu,D)$ is zero. Now the foregoing relations and this observation imply
\begin{equation*}
\mathop{\mathrm{lim\,sup}}_{i\to\infty} \tau_n(\mu_i,D) =  \mathop{\mathrm{lim\,sup}}_{i\to\infty}
\sum_{g\in D'_n} \mu_i(g)\ln\frac{m_{\mu_i}(A^ng)}{\mu_i(g)} \le \sum_{g\in
D'_n}\nu(g)\ln\frac{m(A^ng)}{\nu(g)}-\eps.
\end{equation*}
In view of  \eqref{9,,3} the right-hand part of the latter inequality does not exceed
$\tau_n(\nu,D) -\eps$, which contradicts the continuity of the restriction of~$\tau_n(\,\cdot\,,D)$
to~$M(D'_n)$. \qed \proofskip

Denote by $\Int M(D'_n)$ the set of measures $\mu\in M(D'_n)$ that take positive values on all
elements  $g\in D'_n$. Let $\mu\in \Int M(D'_n)$ and $\nu\in M(D)$. Put
\begin{equation}\label{9,,6}
\tau_n(\nu,\mu,D) := \sum_{g\in D'_n}\nu(g)\ln \frac{m_\mu(A^ng)}{\mu(g)},
\end{equation}
where the measure $m_\mu\in M$ is defined by  equality~\eqref{9,,3}.

\begin{lemma}\label{9..3}
Given a single-valued correspondence\/ $\mu\mapsto m_\mu$ that satisfies\/~\eqref{9,,3}, then for
every measure\/ $\mu_0\in M(D'_n)$ and any number\/ $t>\tau_n(\mu_0,D)$ there exists a
neighborhood\/ $O(\mu_0)$ in the set\/ $M(D)$ and a measure\/ $\mu\in O(\mu_0)\cap\Int M(D'_n)$
such that for all\/ $\nu\in O(\mu_0)$ the following estimate holds:\/ \ $\tau_n(\nu,\mu,D)<t$.
\end{lemma}

\proof. Let $D'_n = \{g_1,\dots,g_k\}$, and  $\mu_1$ be the center of the simplex   $M(D'_n)$. It
is defined by the equalities
$$
 \mu_1(g_i)  = 1/k \quad \textrm{for all}\ \ g_i\in D'_n .
$$
Let us consider the family
$$
 \mu_\theta  = (1-\theta)\mu_0 +\theta\mu_1.
$$
For small positive $\theta$ this family belongs to $\Int M(D'_n)$. Denote by $O_\theta(\mu_0)$ the
set of measures $\nu\in M(D)$ satisfying the inequalities
$$
 \abs{\nu(g_i)-\mu_0(g_i)} <\theta \quad \textrm{for all}\ \ g_i\in D'_n.
$$
Clearly,
$$
\mu_\theta\in O_\theta(\mu_0) .
$$

\medskip\noindent
By Proposition \ref{8..5} the function $\tau_n(\mu_\theta,D)$ depends continuously on the parameter
$\theta\ge 0$ and when $\theta$ is small $\tau_n(\mu_\theta,D)$ differs in a small way from
$\tau_n(\mu_0,D)$. Therefore, it suffice to prove that if  $\nu\in O_\theta(\mu_0)$ then the
difference
\begin{equation}\label{9,,7}
\tau_n(\nu,\mu_\theta,D) -\tau_n(\mu_\theta,D)  = \sum_{g\in D'_n}
\left(\frac{\nu(g)}{\mu_\theta(g)} -1\right) \times
\mu_\theta(g)\ln\frac{m_{\mu_\theta}(A^ng)}{\mu_\theta(g)}
\end{equation}
tends uniformly to zero as $\theta\to 0$. This can be derived in a simple way from
Lemma~\ref{9..2}. Indeed, if  $\mu_0(g)>0$, then for small $\theta$ the corresponding summand in
\eqref{9,,7} has small first factor  and bounded second one; and if  $\mu_0(g) = 0$ then the first
factor is bounded while  the second one is small.\qed

\begin{lemma}\label{9..4}
If\/ $\mu\in\Int M(D'_n)$ and a measure\/ $m_\mu\in M$ satisfies equality\/~\eqref{9,,3}, then
\begin{equation*}
\sum_{g\in D'_n}\frac{\mu(g)}{m_\mu(A^ng)}A^ng\,\le\, 1.
\end{equation*}
\end{lemma}

\proof. Lemma \ref{9..2} implies that $m_\mu(A^ng)>0$ whenever $\mu(g)>0$. Let us take any measure
$m\in M$ and consider the function
\begin{equation*}
\varphi(t)  = \sum_{g\in D'_n}\mu(g)\ln\frac{(1-t)m_\mu(A^ng) + tm(A^ng)}{\mu(g)}, \qquad t\in
[0,1].
\end{equation*}
By \eqref{9,,3} this function attains the maximal value at  $t = 0$. Therefore its derivative at
zero is nonpositive:
\begin{equation*}
\varphi'(0)  = \sum_{g\in D'_n}\mu(g)\frac{m(A^ng) -m_\mu(A^ng)} {m_\mu(A^ng)}  = \sum_{g\in
D'_n}\frac{\mu(g)}{m_\mu(A^ng)}m(A^ng) -1\le 0.
\end{equation*}
Since the measure $m\in M$ is arbitrary this implies the statement of the lemma.  \qed

\begin{lemma}\label{9..5}
If a function\/ $g\in C(X)$ is nonnegative and\/ $A^ng  \equiv  0$, then for every positive
integer\/ $N$ the following equality holds true
$$
A^{N+n}\bigl(e^{S_Ng}\bigr) = A^{N+n}\mathbf 1,
$$
where $S_Ng$ is defined by \eqref{5,,1}.
\end{lemma}

\proof. Evidently, the function $A^{i+n}(g\circ\alpha^i)$ is nonnegative. On the other hand the
homological identity  implies
$$
 A^{i+n} (g\circ\alpha^i)  = A^n(gA^i\mathbf 1) \le A^ng\cdot \norm{A}^i =  0.
$$
Hence $A^{i+n}(g\circ\alpha^i) \equiv 0$ and thus
\begin{equation}\label{9,,8}
A^{N+n}(S_Ng)  = \sum_{i = 1}^N A^{N+n}(g\circ\alpha^{i-1}) \equiv 0 .
\end{equation}
Recall now that the function  $S_Ng$ is nonnegative and bounded. By the Lagrange mean value theorem
we have
$$
 e^{S_Ng}-\mathbf 1 \le e^{\norm{S_Ng}}S_Ng .
$$
Applying the operator $A^{N+n}$ to this inequality and exploiting (\ref{9,,8})we obtain
$$
 A^{N+n}\bigl(e^{S_Ng}\bigr) -A^{N+n}\mathbf 1\le  e^{\norm{S_Ng}}A^{N+n}(S_Ng) =0.
$$
On the other hand, $e^{S_ng}\ge1$, and therefore
$$
 A^{N+n}\bigl(e^{S_ng}\bigr)\ge A^{N+n} \mathbf 1 .  \qed
$$
\proofskip

Now we can prove the Entropy Statistic Theorem itself.

\medskip

\noindent{\bf The proof of the Entropy Statistic Theorem.}
Let us fix a measure $\mu_0\in M$ and an arbitrary number $t>\tau(\mu_0)$. Choose a natural number
$n$ and a partition of unity $D$ in~$C(X)$, satisfying the inequality
$$
 \tau_n(\mu_0,D) <nt.
$$
as above, we denote by $D'_n$  the union of functions $g\in D$ such that $A^ng\not\equiv 0$.

Suppose first that $\mu_0\notin M(D'_n)$. In this case there exists a function $g\in D\setminus
D'_n$ satisfying the inequality $\mu_0(g) >0$. Take the number $\eps = \mu_0(g)/2$ and define the
neighborhood
\begin{equation}\label{9,,9}
O(\mu_0) := \{\,\mu\in M\mid \mu(g) >\eps\,\}.
\end{equation}
Let $x\in X_N(O(\mu_0))$ for some~$N$, where $X_N(O(\mu_0))$ is that defined by (\ref{9,,2}). Then
the empirical measure $\delta_{x,N}$ belongs to $O(\mu_0)$ and hence
$$
 S_Ng(x) = N\delta_{x,N}(g) >N\eps.
$$
Define the function $\chi_N$ by the formula
$$
 \chi_N = e^{C(S_Ng-N\eps)},
$$
where $C$ is arbitrary positive. Clearly, $\chi_N$ majorizes the index function of the set
$X_N(O(\mu_0))$. If in addition $N>n$, then
$$
 \chi_N \le e^{C(S_{N-n}g+n-N\eps)}  =  e^{Cn-CN\eps}e^{S_{N-n} (Cg)} .
$$
This along with Lemma \ref{9..5} implies
\begin{equation}\label{9,,10}
 A^N\chi_N = e^{Cn-CN\eps}A^N\mathbf 1 \le e^{Cn} e^{N(\ln\|A\| -C\eps)}.
\end{equation}
If we choose $C$ so large that  $\ln\norm A -C\eps <t$, then (\ref{9,,10}) implies the Entropy
Statistic Theorem.

To complete the proof it remains  to consider the situation when $\mu_0\in M(D'_n)$. By Lemma
\ref{9..3} in this case there exists a neighborhood $O(\mu_0) \subset M(D)$ and a measure $\mu\in
O(\mu_0)\cap \Int M(D'_n)$ such that for every $\nu\in O(\mu_0)$ we have
$$
 \tau_n(\nu,\mu,D) <nt.
$$

Define the function
\begin{equation}\label{9,,11}
 \psi(x) := \sum_{g\in D'_n}g(x)\ln\frac{\mu(g)}{m_\mu(A^ng)} + \sum_{g\in D\setminus D'_n}
 g(x)\ln 1,
\end{equation}
where the measure $m_\mu\in M$ satisfies equality~\eqref{9,,3}. Note that the second sum in
\eqref{9,,11} is equal to zero. Therefore by comparing \eqref{9,,11} and \eqref{9,,6} we see that
\begin{equation}\label{9,,12}
S_N\psi(x)  = -N\tau_n(\delta_{x,N},\mu,D).
\end{equation}
On the other hand the convexity of the exponent implies
\begin{equation}\label{9,,13}
 e^{\psi(x)}\le \sum_{g\in D'_n}g(x)\frac{\mu(g)}{m_\mu(A^ng)} + \sum_{g\in D\setminus D'_n}g(x).
\end{equation}

Applying the operator $A^n$ to \eqref{9,,13} and exploiting Lemma \ref{9..4} we obtain the estimate
\begin{equation}\label{9,,14}
 A^n\bigl(e^{\psi}\bigr)\le 1.
\end{equation}
Let us introduce the notation
$$
 \psi_k =\psi+\psi\circ\alpha^n+\dots+\psi \circ\alpha^{n(k-1)}.
$$
The homological identity and estimate (\ref{9,,14}) imply
\begin{equation}\label{9,,15}
A^{nk}\bigl(e^{\psi_k}\bigr)  = A^n\bigl(e^{\psi}A^n\bigl( e^{\psi}\dotsc
A^n\bigl(e^{\psi}A^n\bigl(e^{\psi} \bigr)\bigr)...\bigr)\bigr) \le 1.
\end{equation}

By construction the function $\psi(x)$ has a finite norm $\norm\psi$. Consider an integer $N>n$ and
take a natural $k$ such that $N\in [n(k+1),\,n(k+2)]$. Then
\begin{equation}\label{9,,16}
S_N\psi\le S_{nk}\psi +2n\norm\psi  = \sum_{i = 0}^{n-1}\psi_k\circ\alpha^i +2n\norm\psi.
\end{equation}
If $x\in X_N(O(\mu_0))$ then the empirical measure  $\delta_{x,N}$ belongs to $O(\mu_0)$ (see
\eqref{9,,2}) and by the choice of this neighborhood,
$$
 \tau_n(\delta_{x,N},\mu,D)<nt.
$$

\smallskip\noindent In this case  \eqref{9,,12} along with \eqref{9,,16} imply
\begin{equation*}
Nt>\frac{N}{n}\tau_n(\delta_{x,N},\mu,D)  = -\frac{1}{n}S_N\psi(x)\ge -\frac{1}{n}\sum_{i =
0}^{n-1} \psi_k\circ\alpha^i(x) -2\norm\psi
\end{equation*}
for every $x\in X_N(O(\mu_0))$. Therefore, the function $\chi_N$ defined by the formula
\begin{equation*}
\chi_N(x)  := \exp\biggl\{Nt+2\norm\psi +\frac{1}{n}\sum_{i = 0}^{n-1}
\psi_k\circ\alpha^i(x)\biggr\}
\end{equation*}
majorizes the index function of the set $X_N(O(\mu_0))$.

Recall the H\"older inequality: if  $m$ is a linear positive functional on  $C(X)$ and $f_1$,
\dots, $f_n$ are nonnegative functions from $C(X)$ then
$$
 m(f_1\dotsm f_n) \le \bigl(m(f_1^n)\dotsm m(f_n^n)\bigr)^{1/n}.
$$
Applying this inequality to the functional  $m(f)  = \bigl[A^Nf\bigr](x)$ and the function $f =
\chi_N$, we obtain the estimate
\begin{equation}\label{9,,17}
\bigl[A^N\chi_N\bigr](x)\le e^{Nt +2\norm\psi} \prod_{i = 0}^{n-1}
\left[\left(A^Ne^{\psi_k\circ\alpha^i}\right)\!(x)\right]^{1/n}.
\end{equation}
In addition, the homological identity and \eqref{9,,15} imply
\begin{equation}\label{9,,18}
 A^N\Bigl(e^{\psi_k\circ\alpha^i}\Bigr)  = A^{N-i-nk}A^{nk}\bigl(e^{\psi_k}A^i\mathbf
 1\bigr) \le \norm{A}^{N-nk}A^{nk} e^{\psi_k}\le \norm{A}^{N-nk}.
\end{equation}
Estimates (\ref{9,,17}) and (\ref{9,,18}) imply the Entropy Statistic Theorem. \qed\proofskip

Now, by means of this theorem, we can prove the Variational Principle.

\section{The proof of the Variational Principle}\label{10..}

By Theorem \ref{5..5} we have already established that
\begin{equation}\label{10,,1}
\lambda(\varphi)\ge \sup_{\mu\in M_\delta ({\cal C})}\bigl(\mu(\varphi) +\tau(\mu)\bigr),
\end{equation}
and therefore we have to prove the opposite inequality. Its proof needs some auxiliary results.

\begin{prop}\label{10..1}
The sets\/ $M(\cal C)$ and\/ $M_\delta(\cal C)$ are compact in\/ $^*$-weak topology.
\end{prop}

\proof. These sets are closed subsets of the unit ball in the dual space to $\cal C$ and by the
Alaoglu theorem the latter ball is compact. \qed

\begin{prop}\label{10..2}
Let\/ $(X,\alpha)$ be the dynamical system corresponding to a\/ $C^*$-dynamical system\/ $(\cal C,
\delta)$ and\/ $\cal C$ be identified with\/ $C(X)$ \emph{(see\/ \ref{6..1}, \ref{6..4})}. If\/ $U$
is a neighborhood of the set\/ $M_\delta(\cal C)$ in\/ $M(\cal C)$ then there exists a (large)\/
$N$ such that for all\/ $n>N$ and\/ $x\in X$ one has\/ $\delta_{x,n}\in U$, where\/ $\delta_{x,n}$
is the empirical measure defined by\/ \eqref{9,,1}.
\end{prop}

\proof{} is implemented by contradiction. Suppose that there exist sequences $x_k\in X$ and
$n_k\to \infty$ such that $\delta_{x_k,n_k}\notin U$. By the compactness of $M(\cal C)$ the
sequence $\delta_{x_k,n_k}$ possesses a limit point
\begin{equation}\label{10,,2}
 \mu\in M(\cal C)\setminus U.
\end{equation}

\smallskip\noindent
Clearly, for every $f\in \cal C$ the following equality holds
$$
\delta_{x_k,n_k}(f) -\delta_{x_k,n_k}(f\circ\alpha)  = n_k^{-1}(f(x_k) -f(\alpha^{n_k}(x_k)).
$$
As a limit one has $\mu(f)\equiv \mu(f\circ\alpha)$ and so $\mu\in M_\delta(\cal C)$, which
contradicts (\ref{10,,2}).  \qed \proofskip

Now we can finish the proof of the Variational Principle.

Hereafter we exploit the identification $\cal C =C(X)$ explained in \ref{6..1}, \ref{6..4}.

Consider the functional
\begin{equation}\label{10,,3}
 \Lambda(\varphi) := \sup_{\mu\in M_\delta(\cal C)} \bigl(\mu(\varphi) +\tau(\mu)\bigr).
\end{equation}
First we prove the  equality
\begin{equation}\label{10,,4}
\lambda (\varphi) = \Lambda(\varphi),\qquad \varphi \in \cal C .
\end{equation}

\medskip\noindent
Choose arbitrary numbers $c>\Lambda(\varphi)$ and $\eps>0$. On the set of invariant measures
$M_\delta(\cal C)$ we define the functional
$$
 t(\mu) :=  c-\mu(\varphi) .
$$
Clearly,
$$
 t(\mu)> \tau(\mu) .
$$
For every measure $\mu\in M_\delta(\cal C)$ we choose a (small) neighborhood  $O(\mu)$ in $M(\cal
C)$ such that it  satisfies the Entropy Statistic Theorem (with the number $t = t(\mu)$) and, in
addition, for all $\nu\in O(\mu)$ the following inequality holds:
$$
 \nu(\varphi) <\mu(\varphi) +\eps.
$$
Then for every point $x\in X_n(O(\mu))$ we have
\begin{equation}\label{10,,5}
 S_n\varphi(x) = n\delta_{x,n}(\varphi) <n\bigl(\mu(\varphi)+\eps\bigr)  =
n\bigl(c-t(\mu) +\eps\bigr).
\end{equation}
By Proposition \ref{10..1} the set $M_\delta(\cal C)$ is compact. Let us cover it by a family of
neighborhoods of the type mentioned, and let $O(\mu_1),\,\dots,\, O(\mu_k)$ be a finite subcover
consisting of the elements of this family. The Entropy Statistic Theorem assigns to each
neighborhood $O(\mu_i)$ a sequence of functions $\chi_{i,n}$ majorizing the index functions of the
sets $X_n(O(\mu_i))$. Proposition~\ref{10..2} implies that for all sufficiently large  $n$ every
point $x\in X$ lies in a certain  set  $X_n(O(\mu_i))$,\ \ $i = 1,\,\dots,\,k$. Hence,
$$
 \chi_{1,n}+\dotsm+\chi_{k,n}\ge \mathbf 1 .
$$
Now \eqref{10,,5} along with the Entropy Statistic Theorem imply the estimates
\begin{equation*}
A_\varphi^n\mathbf 1  = A^n\bigl(e^{S_n\varphi}\mathbf 1\bigr) \le \sum_{i = 1}^k
A^n\bigl(e^{n(c-t(\mu_i)+\eps)}\chi_{i,n}\bigr) \le \sum_{i = 1}^k e^{n(c-t(\mu_i)+\eps)}
C(t(\mu_i),\mu_i)e^{nt(\mu_i)}.
\end{equation*}
As $n\to\infty$ these estimates imply
$$
 \lambda(\varphi)\le c+\eps.
$$
Since the numbers $c>\Lambda(\varphi)$ and $\eps>0$ are arbitrary the following inequality holds
$$
 \lambda(\varphi) \le \Lambda(\varphi) .
$$
In view of \eqref{10,,1}, \eqref{10,,3} this implies \eqref{10,,4}.

Finally, to complete the proof of (\ref{6,,4}) we recall that by Proposition \ref{8..6} \
$t$-entropy $\tau(\mu)$ is upper semicontinuous on the compact space $M_\delta(\cal C)$. Therefore
the supremum in~\eqref{10,,3} is in fact maximum and thus (\ref{6,,4}) is proved. \qed\proofskip

%
%
%
%
%

\section[Variational Principle for transfer operators with nonnegative weights]{Variational
Principle for transfer operators\\ with nonnegative weights} \label{11..}

In this section we prove a $C^*$-algebraic version of Theorem \ref{5..8}.

Let $({\cal C},\delta)$ be a $C^*$-dynamical system and let $(X, \alpha)$ be the
cor\-res\-pond\-ing dynamical system. Let $A\!:\cal C \to \cal C$ be a \emph{fixed} transfer
operator for  $({\cal C},\delta)$. We  define the family  of operators $A\psi\!:\cal C \to \cal C$,
where $\psi\in \cal C$, as
\begin{equation}\label{11,,1}
A\psi :=A(\psi\,\cdot\,).
\end{equation}
Evidently, if $\psi\ge 0$ then $A\psi$ is a transfer operator.


\begin{rem}\label{11..1}
1) If $\psi>0$ then
\begin{equation*}
 A\psi =A_{\ln\psi},
\end{equation*}
where, as above, $A_\varphi =A(e^\varphi\,\cdot\,)$.

2) Given a transfer operator $A\!:\cal C \to \cal C$ one can also consider the family of operators
$\psi A$ for $\psi\in \cal C$. The homological identity (\ref{6,,3}) implies the equality $\psi A =
A\delta(\psi)$. Thus this family is a subfamily of the one considered in \eqref{11,,1}.
\end{rem}

\begin{theorem}\label{11..2}
{\bf [Variational Principle for transfer operators with nonnegative weights]\,} Let\/ $({\cal
C},\delta)$ be a\/ $C^*$-dynamical system and\/ $(X,\alpha)$ be the corresponding dynamical system.
Let\/ $A\!:\cal C\to \cal C$ be a certain transfer operator for\/ $(\cal C,\delta)$,\ \ $A\psi$ be
a transfer operator defined in\/ \eqref{11,,1}, where\/ $\psi\in \cal C$ and\/ $\psi\ge 0$,
and\/~$\ell(\psi)$ be the logarithm of the spectral radius of\/ $A\psi$. Then the following
variational principle holds true
\begin{equation}\label{11,,2}
 \ell(\psi) =\max_{\mu\in M_\delta(\cal C)} \left(\int_X\ln \psi \,d\mu +\tau(\mu)\right)
\end{equation}
{\rm (where in the notation we exploit the identification of\/ $\cal C$ with the\/
\emph{functional} space\/ $C(X)$ (see\/ \ref{6..1}) and, in particular, we treat\/ $\ln\psi$ as a
function on\/ $X$, and we also identify the functionals\/ $\mu\in M_\delta(\cal C)$ with\/
$\alpha$-invariant Borel probability measures on\/ $X$ (see\/ \ref{6..8})).}
\end{theorem}

\proof. Theorem \ref{6..10} implies \eqref{11,,2} for any positive function $\psi\in{\cal C}=C(X)$.
So we have only to consider the case when $\psi$ has zero points.

Let $\psi_n := \psi + 1/n$. Evidently, the sequence $\psi _n >0$ is monotonically decreasing and
tends in norm to $\psi$. The definition of $\ell(\psi)$  implies that
$$
 \ell(\psi_n)\ge\ell(\psi).
$$
This observation along with upper semicontinuity of the spectral radius implies the existence of
the limit
\begin{equation}\label{11,,3}
 \lim_{n\to\infty}\ell(\psi_n)=\ell(\psi).
\end{equation}

\medskip\noindent
The functions $\psi_n$ satisfy \eqref{11,,2}. Hence for each $\mu\in M_\alpha({\cal C})$ we have
\begin{equation} \label{11,,4}
 \ell(\psi_n)\ge \mu(\ln\psi_n)+ \tau (\mu).
\end{equation}
By the choice of  $\psi_n$ and  Levi's monotonic limit theorem  we have
\begin{equation*}
 \lim_{n\to\infty}\mu(\ln\psi_n)= \mu(\ln\psi),
\end{equation*}
which along with (\ref{11,,3}) and (\ref{11,,4}) means that
\begin{equation}\label{11,,5}
\ell(\psi )\ge \mu(\ln\psi)+\tau(\mu),\qquad \mu\in M_\alpha({\cal C}).
\end{equation}

Now let $\mu_n\in M_\delta ({\cal C})$ be a measure such that
\begin{equation}\label{11,,6}
\ell(\psi_n)= \mu_n(\ln\psi_n)+ \tau(\mu_n),
\end{equation}
and let $\mu_0$ be a limit point of the sequence $\{\mu_n\}$. Evidently, $\mu_0\in M_\delta ({\cal
C})$.

Since $\psi_k >\psi_n$ for $k<n$, equality (\ref{11,,6}) implies also that
\begin{equation*}
 \ell(\psi_n)\le \mu_n(\ln\psi_k)+ \tau(\mu_n), \qquad k<n.
\end{equation*}
This inequality along with (\ref{11,,3}) and the upper semicontinuity of $\tau(\mu)$ (Proposition
\ref{8..5}) implies
\begin{equation*}
\ell(\psi)\le \mu_0(\ln\psi_k)+ \tau (\mu_0).
\end{equation*}

\medskip\noindent
Finally, passing to the limit and applying Levi's monotonic limit theorem, we get
\begin{equation*}
 \ell(\psi)\le \mu_0(\ln\psi)+ \tau (\mu_0).
\end{equation*}
Which along with (\ref{11,,5}) proves (\ref{11,,2}). \qed

\section[Additional properties and  alternative definition of $\boldsymbol t$-entropy]{Additional
 properties and  alternative\\ definition of $\boldsymbol t$-entropy} \label{12..}

The Variational Principle tells us in particular that $t$-entropy plays the principal  role in the
spectral analysis of transfer operator and this section is devoted to the description of a number
of additional properties of $t$-entropy that have not been exploited in the proof of the
Variational Principle. Along with this we will also present a certain alternative definition of
$t$-entropy.

\begin{prop}\label{12..1}
If\/ $\mu\in M_\delta(\cal C)$ then\/
\begin{equation}\label{12,,1}
 \tau_{n+k} (\mu)\le \tau_n(\mu) +\tau_k(\mu).
\end{equation}
Therefore the\/ $t$-entropy of a\/ $\delta$-invariant measure\/ $\mu$ can be defined as the limit
$$
\tau(\mu) =\lim_{n\to\infty} \frac{\tau_n(\mu)}{n}.
$$
\end{prop}

\proof. Let us consider two partitions of unity  $D,\,E$ in the algebra~$\cal C$. For every two
elements  $g\in D$ and $h\in E$ we define the function
$$
 f_{gh} :=g\,h\circ\alpha^k .
$$
Consider the new partition of unity
$$
 C:=\{\,f_{gh} \mid g\in D,\ h\in E\,\} ,
$$
and introduce the notation
$$
 D_\mu :=\{\,g\in D\mid \mu(g)>0\,\} ,
$$
and
$$
 E_\mu :=\{\,h\in E\mid \mu(h)>0\,\} .
$$

\medskip

Suppose first that $A^nh = 0$ for a certain $h\in E_\mu$. Then by the homological identity we have
$$
 A^{n+k}f_{gh} =A^n(hA^kg) = 0.
$$
Therefore,
$$
 \tau_n(\mu,E) =\tau_{n+k}(\mu,C) = -\infty
$$

\smallskip\noindent
 and the proposition is proved.

Now let all the elements $A^nh$, where  $h\in E_\mu$, be nonzero. Then there exist measures  $m\in
M(\cal C)$, satisfying the condition  $m(A^nh)>0$ for all  $h\in E_\mu$. These measures by the
$\delta$-invariance of the measure  $\mu$ and concaveness of the logarithm function satisfy the
relations
\begin{gather*}
\sum_{f_{gh}\in C}\mu(f_{gh})\ln\frac{m(A^{n+k}f_{gh})}
{\mu(f_{gh})} -\sum_{h\in E_\mu}\mu(h)\ln\frac{m(A^nh)}{\mu(h)}\displaybreak[0]\\[3pt]
 =\sum_{g\in D_\mu}\mu(g)\sum_{h\in E_\mu}\frac{\mu(f_{gh})}{\mu(g)}\ln
\frac{m(A^{n+k}f_{gh})\mu(h)}{\mu(f_{gh})m(A^nh)}\displaybreak[0] \\[3pt]
\le \sum_{g\in D_\mu}\mu(g)\ln\biggl(\frac{1}{\mu(g)}
\sum_{h\in E_\mu}\frac{m(A^{n+k}f_{gh})\mu(h)}{m(A^nh)}\biggr)\displaybreak[0] \\[3pt]
\qquad =\sum_{g\in D_\mu}\mu(g)\ln\frac{m'(A^kg)}{\mu(g)}, \qquad \text{where}\quad m'(\,\cdot\,)
=\sum_{h\in E_\mu}\frac{m(A^n(h\,\cdot\,)) \mu(h)}{m(A^nh)}.
\end{gather*}
One can easily verify that $m'\in M(\cal C)$. By varying the measure $m$ in the latter relation we
obtain the inequality
$$
 \tau_{n+k}(\mu,C) \le \tau_n(\mu,E) +\tau_k(\mu,D) ,
$$
which implies (\ref{12,,1}). \qed

\begin{I}\label{12..2}
For any partition of unity $D$ in the algebra $C(X)$ and any positive number $\eps$ we denote by
$W(D,\eps)$ the set of all partitions of unity  $E$ in  $C(X)$, satisfying the following condition:
if $g\in D$ and some function $h\in E$ does not vanish at points $x_1,x_2\in X$ then $\abs{g(x_1) -
g(x_2)}<\eps$. In other words, each function $g\in D$ has small oscillation on the support of any
$h\in E$.
\end{I}

\begin{prop}\label{12..3}
Given a measure\/ $\mu\in M(C(X))$ and a partition of unity\/ $D$ in\/ $C(X)$, and a real number\/
$t>\tau_n(\mu,D)$, there exists a (small)\/ $\eps>0$ such that for every partition of unity\/ $E\in
W(D,\eps)$ the estimate\/ $\tau_n(\mu,E) <t$ holds true.
\end{prop}

\proof. Let
$$
 D_\mu :=\{\,g\in D\mid \mu(g)>0\,\}.
$$

\medskip\noindent
Take a (small) positive number $\eps$ such that for every $g_0\in D_\mu$ the following inequality
holds
\begin{equation}\label{12,,2}
\mu(g_0)\ln\frac{\sqrt\eps +\eps\norm{A}^n}{\mu(g_0)}\, + \sum_{g\in D_\mu\setminus\{g_0\}}
\mu(g)\ln \frac{2\norm{A}^n}{\mu(g)} \,<\,t.
\end{equation}

\medskip\noindent
Consider any partition of unity $E\in W(D,\eps)$, where $W(D,\eps)$ is defined in \ref{12..2}. Let
$$
 E_\mu :=\{\,h\in E\mid \mu(h)>0\,\} .
$$
Then for every measure  $m\in M(C(X))$ the concavity of the logarithm implies that
\begin{align}\notag
\sum_{h\in E_\mu}\mu(h)\ln\frac{m(A^nh)}{\mu(h)} &=\sum_{g\in D_\mu}\mu(g)\sum_{h\in
E_\mu}\frac{\mu(gh)}{\mu(g)}\ln
\frac{m(A^nh)}{\mu(h)}\\[3pt]\label{12,,3}
& \le \sum_{g\in D_\mu}\mu(g)\ln\Biggl(\frac{1}{\mu(g)} \sum_{h\in
E_\mu}\frac{\mu(gh)}{\mu(h)}m(A^nh)\Biggr).
\end{align}
If $g\in D_\mu$ and $h\in E_\mu$, then the definition of $W(D,\eps)$ implies
$$
 \frac{\mu(gh)}{\mu(h)}h\le (g+\eps)h.
$$
Substituting this inequality into \eqref{12,,3} we obtain
\begin{equation}\label{12,,4}
\sum_{h\in E_\mu}\mu(h)\ln\frac{m(A^nh)}{\mu(h)} \le \sum_{g\in
D_\mu}\mu(g)\ln\frac{m(A^n(g+\eps))}{\mu(g)}.
\end{equation}

Now it suffice to check that for small $\eps$ the sum in the right-hand part of~\eqref{12,,4} does
not exceed  $t$. Indeed, if it is so then passing to the supremum over $m\in M(C(X))$ we obtain the
estimate $\tau_n(\mu,E)\le t$ we are looking for.

If there exists a function $g_0\in D_\mu$ such that $m(A^ng_0) <\sqrt\eps$, then evidently
$$
 m\bigl(A^n(g_0+\eps)\bigr) <\sqrt\eps+ \eps\norm{A}^n,
$$
and the whole of the sum \eqref{12,,4} does not exceed $t$ in view of \eqref{12,,2}.

In the opposite case for every function  $g\in D_\mu$ the following estimate holds
$$
m\bigl(A^n(g+\eps)\bigr) \le m(A^ng)\bigl(1+\sqrt\eps \norm{A}^n \bigr) .
$$
It implies that the sum in the right-hand part of \eqref{12,,4} does not exceed
$$
\tau_n(\mu,D) +\ln\bigl(1+\sqrt\eps\norm{A}^n\bigr).
$$
And the latter expression is less than $t$ whenever $\eps$ is small enough. \qed \proofskip

Recalling the identification of $\cal C$ with $C(X)$ we are giving now an alternative definition of
$t$-entropy. Let $({\cal C},\delta)$ be a $C^*$-dynamical system and $(X, \alpha)$ be the
corresponding dynamical system. Let~$A$ be a transfer operator in $C(X)$ (we identify $\cal C$ with
$C(X)$)
 and let $A^*\!:C^*(X)\to C^*(X)$ be the adjoint operator to~$A$. The
formulae~\eqref{5,,8},~\eqref{5,,9} can be rewritten in the following way:
\begin{gather}\label{12,,5}
\tau(\mu) :=\inf_{n}\frac{\tau_n(\mu)}{n}\,,\qquad
\tau_n(\mu) :=\inf_D\tau_n(\mu,D),\\[6pt]
\tau_n(\mu,D) :=\sup_{m\in M(C(X))}\sum_{g\in D}\mu(g) \ln\frac{A^{*n}m(g)}{\mu(g)}\,.\label{12,,6}
\end{gather}
Here  $D$ denotes a \emph{continuous} partition of unity on~$X$.

In fact these formulae have sense for arbitrary \emph{measurable} partitions of unity as well.

Let us consider any measure $\mu\in M(C(X))$. A finite Borel partition $\{G_1,\dots,G_k\}$ of the
\emph{space}~$X$ will be called $\mu$\emph{-proper}, if the boundary of each set $G_i$ has zero
measure~$\mu$.

A Borel partition of \emph{unity} $D=\{g_1,\dots,g_k\}$ on $X$ will be called $\mu$\emph{-proper}
if it consists of index functions of sets $G_i\subset X$ forming a $\mu$-proper partition of~$X$.

\begin{prop}\label{12..5}
Given an open cover\/ $W$ of a Hausdorff compact space\/~$X$ and a measure\/ $\mu\in M(C(X))$,
there exists a finite\/ $\mu$-proper partition of\/~$X$ refined in\/~$W$.
\end{prop}

\proof. Without loss of generality one can assume that the open cover~$W$ is finite. Let us
consider the case when it consists of two sets $U_0$, $U_1$. Suppose that the closed sets
$\overline U_i =X\setminus U_i$ are nonempty  (in the opposite case the statement is trivial). By
the Urysohn's Lemma there exists  a continuous function  $f\!:X\to [0,1]$ vanishing on $\overline
U_0$ and equal to unit on $\overline U_1$. By the finiteness of the measure  $\mu$ there exists a
number $c\in (0,1)$ such that $\mu(f^{-1}(c)) =0$. Now one can take the partition of the space~$X$
consisting of two sets  $G_0 =f^{-1}([0,c))$ and $G_1 =f^{-1}([c,1])$. In the general case one
should apply the induction by the number of elements of the cover and use on every step the
construction described. \qed \proofskip


\begin{prop}\label{12..6}
The values of the functionals\/ $\tau(\mu)$ and\/ $\tau_n(\mu)$, where\/ $\mu\in M(C(X))$, do not
change if in formulae\/~\eqref{12,,5} and\/~\eqref{12,,6} one uses\/ $\mu$-proper partitions of
unity instead of continuous partitions of unity\/~$D$.
\end{prop}

\proof. Note first that Proposition~\ref{12..3} holds true if in this proposition one considers
arbitrary Borel partitions of unity instead of continuous ones and defines the functional
$\tau_n(\mu,D)$ by means of formula \eqref{12,,6}. Its proof stays the same, one should simply
replace everywhere $m(A^n\,\cdot\,)$ by~$A^{*n}m(\,\cdot\,)$.

Define the functional
$$
 \tau'_n(\mu) :=\inf_D \tau_n(\mu,D),
$$
where the infimum is taken over all $\mu$-proper partitions of unity~$D$. Propositions \ref{12..5}
and \ref{12..3} (for the Borel partitions) imply
$$
\tau'_n(\mu)\le \tau_n(\mu) .
$$
So it suffice to prove the opposite inequality.

Fix a number $\eps>0$. Choose a $\mu$-proper partition of unity $E =\{h_1,\dots,h_k\}$ such that
$$
\tau_n(\mu,E) <\tau'_n(\mu) +\eps .
$$
By definition the sets $V_i =h_i^{-1}(1)$ form a $\mu$-proper partition of $X$.
 Denote by
$\Gamma$ the union of the boundaries of these sets. Clearly, $\mu(\Gamma) =0$. Let~$\delta_\eps$ be
a positive and small enough number (it will be indicated below how small it should be). By the
regularity of the measure $\mu$ the set  $\Gamma$ possesses a neighborhood~$O(\Gamma)$ such that
$$
 \mu(O(\Gamma)) <\delta_\eps.
$$
Let us consider the open cover of  $X$ consisting of the sets
$$
 O(\Gamma),\,V_1\setminus \Gamma,\, \dots,\,V_k\setminus \Gamma
$$
and choose  a continuous partition of unity  $D =\{g_0,g_1, \dots, g_k\}$ refined in this partition
(in such a way that  $\mathop{\mathrm{supp}} g_0\subset O(\Gamma)$ and $\mathop{\mathrm{supp}} g_i
\subset V_i\setminus\Gamma$). Then
$$
 g_i\le h_i\le g_i +g_0, \qquad i=1,\,\dots,\,k .
$$
In addition,
$$
 \mu(g_0)\le \mu(O(\Gamma)) <\delta_\eps.
$$
Therefore,
\begin{equation}\label{12,,7}
\bigl|\mu(g_i) -\mu(h_i)\bigr| <\delta_\eps, \qquad i=1,\,\dots,\,k.
\end{equation}

\medskip\noindent
Let $E'_n$ be the set of all the functions $h_i\in E$ for which  there exists a measure $m'\in
M(C(X))$ such that $A^{*n}m'(h_i) >0$. Fix a measure  $m'\in M(C(X))$ satisfying the condition
$$
 A^{*n}m'(h_i) >0 \quad \textrm{for  all}\ \ h_i\in E'_n.
$$

Let $m\in M(C(X))$. Consider the measure
$$
 m_\eps =(1-\eps)m+\eps m'.
$$
It satisfies the inequalities
\begin{gather}\nonumber
\ln(1-\eps) +\sum_{i=0}^k\mu(g_i)\ln\frac{m(A^ng_i)}{\mu(g_i)}
\le\sum_{i=0}^k\mu(g_i)\ln\frac{m_\eps(A^ng_i)}{\mu(g_i)} \\[3pt]
 \le\mu(g_0)\ln\frac{\norm{A}^n}{\mu(g_0)} +\sum_{i=1}^k
\mu(g_i) \ln\frac{A^{*n}m_\eps(h_i)}{\mu(g_i)}.  \label{12,,8}
\end{gather}

It suffice to check that the whole of the sum \eqref{12,,8} does not exceed~$\tau_n(\mu,E) +\eps$.
Since, if it is so, then taking the supremum over $m\in M(C(X))$ one obtains the estimate
$$
 \ln(1-\eps) +\tau_n(\mu,D)\le \tau_n(\mu,E) +\eps <\tau'_n(\mu) +2\eps,
$$
and by the arbitrariness of  $\eps$ it implies that
$$
 \tau_n(\mu)\le \tau'_n(\mu).
$$

The first summand in \eqref{12,,8} depends continuously on the value of $\mu(g_0)$. Since $\mu(g_0)
<\mu(O(\Gamma)) <\delta_\eps$, this summand can be made arbitrarily small by the choice of
$\delta_\eps$. Definition~\eqref{12,,6} implies that if one replaces in the second summand
in~\eqref{12,,8} all the functions $g_i$ by $h_i$, then he obtains the expression not
exceeding~$\tau_n(\mu,E)$. Therefore it suffice to check that the second summand in \eqref{12,,8}
changes in a small way under this replacement.

Let us consider each summand in the right-hand sum of~\eqref{12,,8} separately.

If  $h_i\in E'_n$ then by the construction one has
$$
 A^{*n}m_\eps(h_i) \ge \eps A^{*n}m'(h_i) >0,
$$
and on the other hand
$$
 A^{*n}m_\eps(h_i) \le\norm{A}^n.
$$

\medskip\noindent
These two inequalities imply that the corresponding summand in \eqref{12,,8} depends continuously
on the value $\mu(g_i)$ and therefore will change in a small way on the replacement of $\mu(g_i)$
by $\mu(h_i)$ (in view of~\eqref{12,,7}).

If  $h_i\notin E'_n$ and concurrently  $\mu(h_i)>0$, than one has $\,\mu(g_i) >0$ (again in view
of~\eqref{12,,7}). So in this case the corresponding summand in \eqref{12,,8} and the whole of the
sum \eqref{12,,8} are equal to~$-\infty$.

Finally, if $h_i\notin E'_n$ and $\mu(h_i) =0$, then $\mu(g_i)\le \mu(h_i) =0$ and the
corresponding summand in~\eqref{12,,8} is equal to zero. As a result we have that the whole of the
sum~\eqref{12,,8} does not exceed $\tau_n(\mu,E) +\eps$ as soon as the number~$\delta_\eps$
in~\eqref{12,,7} is sufficiently small. \qed

\medskip

This proposition naturally leads to the following

\begin{defn}\label{proper}

Let $({\cal C},\delta)$ be a $C^*$-dynamical system and $(X, \alpha)$ be the corresponding
dynamical system. Let~$A$ be a transfer operator in $C(X)$ (we identify $\cal C$ with $C(X)$) and
let $A^*\!:C^*(X)\to C^*(X)$ be the adjoint operator to~$A$. \ \ \emph{$T$-entropy} is the
functional $\tau$ on the set $M_\alpha$ (of all $\alpha$-invariant Borel probability measures on
$X$) defined by the following formulae
\begin{gather}\label{12,,5+}
\tau(\mu) :=\inf_{n}\frac{\tau_n(\mu)}{n}\,,\qquad
\tau_n(\mu) :=\inf_D\tau_n(\mu,D),\\[6pt]
\tau_n(\mu,D) :=\sup_{m\in M(C(X))}\sum_{g\in D}\mu(g) \ln\frac{A^{*n}m(g)}{\mu(g)}\,.\label{12,,6+}
\end{gather}
Here  $D$ denote $\mu$-proper partition of unity.
\end{defn}

Proposition \ref{12..6} tells us that this definition is equivalent to  definition~\eqref{6,,6},
\eqref{6,,7}.

\begin{rem}
To finish this section let us note that the first to be introduced was the definition of
$t$-entropy leaning on  $\mu$-proper partitions, and namely this definition is presented
in~\cite{ABL-2000-1,ABL-2000-2,ABL-2001,ABL-2005}.
\end{rem}

\section[Weighted shift operators and Variational Principle]{Weighted shift operators and\\
Variational Principle} \label{13..}

Example \ref{7..6} shows the tight interrelation between transfer operators and weighted shift
operators in $L^1 (Y,m)$. Developing the idea of this example we go further and present in this
section the Variational Principle for the spectral radius of weighted shift operators acting in
$L^p(Y,m)$ spaces.

\medskip

\noindent
\textbf{Model example.} Just as in Example \ref{7..6}, let $(Y,\goth A,m)$ be a measurable space
with a $\sigma$-finite measure $m$, and $\beta$ be a measurable mapping of $Y$ into itself
satisfying the condition
\begin{equation}\label{13,,1}
m\bigl(\beta^{-1}(G)\bigr) \le Cm(G), \qquad G\in\goth A,
\end{equation}

\medskip\noindent
where the constant $C$ does not depend on $G$. Let us consider the space $L^p(Y,m)$, where $1\le p
\le \infty$. Set the \emph{shift operator} $T$  (generated by the mapping $\beta$) by the formula
\begin{equation}\label{13,,2}
[Tf](x) = f\bigl(\beta(x)\bigr), \qquad f\in L^p(Y,m).
\end{equation}
Inequality \eqref{13,,1} implies that the norm of this operator does not exceed~$C^{1/p}$.

Note that in the case when  $p=\infty$ this operator defines an \emph{endomorphism} of the algebra
$L^\infty (Y,m)$, that is
$$
 T(fg) =Tf\cdot Tg, \qquad f,g\in L^\infty(Y,m).
$$

Just as in Example  \ref{7..6} we define the linear operator $A\!: L^\infty(Y,m) \to L^\infty(Y,m)$
by means of the identity
\begin{equation}\label{13,,3}
 \int_Y f\cdot g\circ\beta\,dm\,\equiv\,\int_Y(Af)g\,dm,\qquad g\in L^1(Y,m),\ \ f\in L^\infty(Y,m)
\end{equation}
(in other words, $A$ is adjoint to the shift operator $T$ on $L^1(Y,m)$).

The definition implies that  $A$ is a positive operator and it satisfies
 the homological identity
\begin{equation}
\label{13,,4}
 A\bigl((Tf)g\bigr) =fAg, \qquad f,g\in L^\infty(X,m).
\end{equation}

\medskip\noindent
Therefore $A$ is a transfer operator (for the $C^*$-dynamical system $(L^\infty(Y,m),T)$). And in
the case of a $\beta$-invariant measure $m$ it is a conditional expectation operator.

Recall (see Section \ref{6..}) that the $C^*$-dynamical system $(L^\infty(Y,m),T)$ can be
canonically identified with the common dynamical system  $(X,\alpha)$, where the compact space~$X$
is the maximal ideal space of the algebra $L^\infty(Y,m)$, and the mapping  $\alpha\!:X\to X$ is
continuous. Under this identification the functions from $L^\infty(Y,m)$ are identified with the
elements of $C(X)$, and the shift mapping  $f\mapsto f\circ\beta$ on  $L^\infty (Y,m)$ is
identified with the shift mapping $g\mapsto g\circ\alpha$ on $C(X)$. Finally, as it was noted in
Remark \ref{7..4}, the set  $M_\alpha(X)$ of all \hbox{$\alpha$-in}\-variant probability measures
on $X$ is identified with the set of all $\beta$-invariant finitely additive probability measures
on $Y$ which are absolutely continuous with respect to $m$. We will denote the latter set by
$M_\beta(Y,m)$.

Since the set $M_\beta(Y,m)$ consists of \emph{finitely} additive measures one can come across
certain difficulties when defining the integrals by these measures for unbounded functions, and
namely such integrals are needed in the next theorem. Fortunately, we can introduce them in a
rather natural way by using the corresponding measures on $C(X)$. Namely, let $\psi\in
L^\infty(Y,m)$ and $\hat\psi\in C(X)$ be its Gelfand transform, let also $\mu\in M_\beta(Y,m)$ and
$\hat\mu\in M_\alpha(X)$ be the corresponding measure mentioned above. Then we set
\begin{equation}\label{13,,5}
 \int_{Y}\ln\abs{\psi}d\mu := \int_{X}\ln\bigl|\hat\psi\bigr|\,d\hat\mu.
\end{equation}

For any $\psi\in L^\infty(Y,m)$ the operator $\psi T$ acting on $L^p(Y,m)$ and given by
$$
 \psi T\!: f\mapsto \psi\cdot Tf
$$
will be called a \emph{weighted shift} operator (with the weight $\psi$). Note, in particular, that
\begin{equation}\label{13,,6}
 T\psi = T(\psi)T, \qquad \psi\in L^\infty(Y,m).
\end{equation}

\begin{theorem}\label{13..1}{\bf [Variational principle for weighted shift operators]\,}
For the  spectral radius of the operator\, $\psi T\!:L^p(Y,m)\to L^p(Y,m)$,\ \ $1\le p<\infty$, the
following variational principle holds:
\begin{equation}\label{13,,7}
 \ln r(\psi T) =\max_{\mu\in M_\beta(Y,m)}\left(\int_{Y}\ln|\psi|\,d\mu +\frac{\tau(\mu)}{p}\right),
\end{equation}
where $\tau(\mu)$ is the $t$-entropy assigned to the transfer operator\/ \eqref{13,,3} and the
integral is understood in the sense of\/ \eqref{13,,5} .
\end{theorem}

\proof. In the case when $p=1$ the operator $A$ is adjoint to $T$. Therefore, the operator $A\psi$
is adjoint to $\psi T$, and they have the same spectral radii. Thus, Theorem~\ref{11..2} implies
\begin{equation}\label{13,,8}
\ln r(\psi T) =\ln r(A\psi) =\max_{\mu\in M_\beta(Y,m)}\left(\int_{Y}\ln|\psi|\,d\mu
+\tau(\mu)\!\right).
\end{equation}

Now let us consider the case  $p>1$. Note that for every function  $f\in L^p(Y,m)$ one has
\begin{equation*}
\int_{Y}\bigl|(\psi T)^n f\bigr|^p\,dm =\int_Y \prod_{i=0}^{n-1}\abs{\psi\circ\beta^i}^p
\bigl|f\circ\beta^n\bigr|^p\,dm =\int_Y\bigl(|\psi|^p T\bigr)^n |f|^p\,dm.
\end{equation*}
So
\begin{equation}\label{13,,9}
\bigl\|(\psi T)^n\bigr\|_{L^p(Y,m)}^p =\bigl\|(|\psi|^p T)^n\bigr\|_{L^1(Y,m)}
=\bigl\|\left(A|\psi|^p\right)^n\bigr\|_{L^\infty(Y,m)},
\end{equation}

\medskip\noindent
and therefore
$$
 p\ln r(\psi T) =\ln r\bigl(|\psi|^p T\bigr),
$$

\medskip\noindent
where the operator $\psi T$ acts on the space $L^p(Y,m)$, and the operator $|\psi|^pT$ acts on he
space  $L^1(Y,m)$. Substituting in \eqref{13,,8} the function $|\psi|^p$ in place of $\psi$ and
dividing the result by $p$ we obtain \eqref{13,,7}. \qed

\begin{rem}\label{13..2}
If $p\to \infty$ then formula (\ref{13,,7}) transforms into the formula
$$
 \ln r(\psi T) =\max_{\mu\in M_\beta(Y,m)} \left(\int_{Y}\ln|\psi|\,d\mu \right).
$$
This restores the Variational Principle for the space $L^\infty (Y,m)$ (cf.~Introduction).
\end{rem}

\noindent{\bf Abstract weighted shifts. Axiomatization.} The study of the weighted shift operators
naturally needs the usage of complex Banach algebras rather than the real ones. Note that the
statement of Theorem~\ref{6..2} is valid not only for the base algebra but also for any semisimple
commutative Banach algebra with an identity (see \cite{Serin}) (we recall that a commutative Banach
algebra is called semisimple if it has zero radical, that is the intersection of all maximal ideals
is  zero; in this case the Gelfand transform is an isomorphism).
 Consideration of the weighted shift operators in  $L^p(Y,m)$ spaces presented above
and, in particular, their interrelation with  the naturally arising transfer operators (see
(\ref{13,,4}) and (\ref{13,,9})) and (\ref{13,,6}) along with the mentioned description of
endomorphisms of semisimple commutative Banach algebras makes it natural the introduction of the
following  axiomatization  of  the weighted shift operators acting in $L^p$ type spaces.

\begin{defn}\label{13..3}
Let ${\cal B}\subset L(B)$ be a semisimple commutative subalgebra of the algebra $L(B)$ of all
linear continuous operators acting on a Banach space $B$ and containing the identity operator
$\mathbf 1$; and let $\delta$ be an endomorphism of $\cal B$ such that $\delta(\mathbf 1) =\mathbf
1$ (and hence by Theorem \ref{6..2} having the form \eqref{6,,1}). Let ${\cal C}$ be a certain
functional base algebra on $X$ (where $X$ is the maximal ideal space of ${\cal B}$) containing
functions of the form $\mathop{\mathrm{Re}}\varphi$ for all $\varphi\in\cal B$, and such that
$({\cal C},\delta)$ is a $C^*$-dynamical system. Let also $A\!:\cal C\to\cal C$ be a certain
transfer operator for $({\cal C}, \delta)$.

We will say that an operator $T\in L(B)$ is an  \emph{abstract shift operator} (associated with
$\delta$ and $A$) and $\psi T$,\ \ $\psi\in {\cal B}$, is an  \emph{abstract weighted shift
operator} (in a space of\/ $L^p$ type,\ \ $1\le p<\infty$), if

a) the equality \/ $T \varphi = \delta(\varphi )T$,\ \ $\varphi\in{\cal B}$ holds and

b) it holds the identity
\begin{equation}\label{13,,10}
 \bigl\|(\varphi T)^n\bigr\|_{L(B)}  = \bigl\|\left(A|\varphi|^p\right)^n\!
 \mathbf 1\bigr\|^{1/p}_{\cal C}, \qquad \varphi\in \cal B.
\end{equation}
\end{defn}

\begin{rem}\label{13..4}
1) Since $\cal C$ is a selfadjoint part of a $C^*$-algebra it follows that if $\vert\varphi\vert
\in \cal C$ then $\vert\varphi\vert^p\in \cal C$ for every $1\le  p<\infty$.

2) The model example presented above is a special case of a general scheme and therewith the
\emph{complexity} of the spectral radius calculation in the general scheme is \emph{equal} to that
of a model example  (cf.~(\ref{13,,10}) and (\ref{13,,9})).

3) Recalling Remarks \ref{7..4} \, 4) \, and \,  \ref{7..7} \, 3) \, we have to stress that
 in general given a concrete functional algebra and its  endomorphism then a transfer operator is \emph{not}
necessarily associated with \emph{a point mapping} of the domain and therefore abstract shift operator and abstract weighted shift operator
do not have to originate from any mapping of the domain.

4) Observe that for ${\cal C} = C(X)$ and  \emph{any} positive left inverse $A$ to $\delta$ (that
is $A$ is a conditional expectation operator) there exists a realization of the objects mentioned
in Definition \ref{13..3} (for $p=1$). Indeed, let $E = C(X)^*$,\ \ $T= {A}^*$ and for any $\psi\in
C(X)$ we define the operator $\psi\!: C(X)^* \to C(X)^*$ by the formula
$$
 (\psi\xi)f = \xi(\psi f),\qquad \xi\in C(X)^*,\ \ f \in C(X),
$$
where $[\psi f](x) = \psi(x)f(x)$ in the right-hand part.

Routine check shows that for $\psi$ and $T$ defined in this way all the conditions of
Definition~\ref{13..3} (for $p=1$) are satisfied.

5) If  $A$ is an invertible conditional expectation operator in $\cal C$ then $\Vert Af\Vert =\Vert
f\Vert$,\ \ $f\in \cal C$ (recall the reasoning in the proof of Proposition \ref{8..3}) and
therefore formula  (\ref{13,,10}) transforms  (for any  $p$) into the formula
\begin{equation}\label{13,,11}
\bigl\Vert(\psi  T)^n\bigr\Vert_{L(B)} = \left\Vert\prod_{k=0}^{n-1}\vert\psi \vert\circ\alpha^k
\right\Vert_{\cal C}.
\end{equation}
Precisely according to this formula there was calculated the norm $\Vert(\psi T)^n\Vert$  in the
process of deducing the variational principle (\ref{1,,4}).
\end{rem}

\begin{rem}\label{13..5}
We would like to emphasize that the construction of an appropriate transfer  operator $A$ in the
model example, namely, the operator by means of which one can calculate the norm of $(\psi T)^n$
with the help of formula (\ref{13,,9}) shows that the operator required should contain information
as on the initial measure $m$ so also on the interrelation between this measure and the mapping
$\beta$. Therefore the choice of $A$ in Definition \ref{13..3} reflects in essence an abstract way
of recording  the corresponding information.

To clarify this remark we present the following observation.

Let $T$ be the shift operator on $L^p(Y, m)$ defined in Model Example. And let $\xi$ be the
partition of $Y$ formed by the inverse images of $\beta$, that is
$$
 \xi = \{\beta^{-1}(y)\}_{y\in Y}.
$$
We denote by  $\xi(x)$ the element of $\xi$ containing $x$. Consider the canonical factor space
$(Y_{\xi},\goth A _{\xi},m_{\xi})$ corresponding to the partition $\xi$ and the set of canonical
conditional measures $m^{\tau}(y)$, where $\tau =\xi(x)$ for some $x$. The measures $m^{\tau}$ are
probability measures (that is $m^{\tau}(\tau) =1$ for each $\tau$) and are defined by the equality
$$
 \int_Y f(x)\,dm(x) =\int_{Y_{\xi}}dm_{\xi}(\tau)\int_{\tau}f(y)\,dm^{\tau}(y),\qquad f\in L^1(Y,m)
$$
(the details see, for example, in \cite{MarIngl}, 1.5.)

Define the conditional expectation operator $E$ in the space $ L^{\infty}(Y, m )$ by the formula
\begin{equation}\label{13,,12}
 [E\psi](x) = \int_{\beta^{-1}(x)}\psi(y)\,dm^{\beta^{-1}(x)}(y).
\end{equation}

It is clear that if $0\le\psi\in L^{\infty}(Y,m)$ then ${E}(\psi)\ge 0$ and $E$ satisfies the
homological identity
\begin{equation}\label{13,,13}
 E\bigl((\varphi\circ\beta)\psi\bigr) =\varphi E\psi \quad \textrm{for all}\ \
 \varphi,\psi\in L^{\infty}(Y,m).
\end{equation}

\medskip\noindent
One may note also that ${E}{\mathbf 1} = \mathbf 1$, which implies (in view of (\ref{13,,13})) that
$E$ is a conditional expectation operator (in the sense of Definition \ref{6..5}, where as the base
algebra~$\cal C$ we take the algebra of real-valued functions in $L^{\infty}(Y m)$).

Let $\beta(m)$ be the measure defined by the equality
$$
 \beta(m)(G):= m(\beta^{-1}(G)), \qquad G\in \goth A,
$$ and let $d\beta(m)/dm$ be the Radon--Nicodim derivative of $\beta (m)$ with respect to $m$.
Inequality (\ref{13,,1}) implies that $d\beta(m)/dm\in L^{\infty}(Y,m)$ and its norm does not
exceed $C$. Clearly, the operator $A$ in (\ref{13,,3}) satisfies the equality
\begin{equation}\label{13,,14}
 A=\frac{d\beta (m)}{dm}E,
\end{equation}
which describes the subtle interrelation between $A$, $m$ and $\beta$.
\end{rem}

The Variational Principle for transfer operators  obtained in the foregoing sections and the
reasoning in the proof of Theorem \ref{13..1} lead to the next variational principle for the
abstract weighted shift operators.

\begin{theorem}\label{13..6}
{\bf [Variational Principle for the abstract weighted shift operators in \boldmath{$L^p$} type
spaces]\,} Let\/ $\psi T$ be an operator satisfying equality\/ \eqref{13,,10} (in particular, $\psi
T$ can be an abstract weighted shift operator described in Definition \ref{13..3}). For the
spectral radius of\/ $\psi T$ the following variational principle  holds
\begin{equation}\label{13,,15}
\ln r(\psi T) =\max_{\mu\in M_\delta(\cal C)} \left(\int_{{X}}\ln\vert\psi\vert\,d\mu
+\frac{\tau(\mu)}{p}\right),
\end{equation}
where\/ ${X}$ is the maximal ideal space of\/ $\cal B$.
\end{theorem}

\proof. Recalling (\ref{11,,1})  and (\ref{13,,10}) and applying Theorem \ref{11..2} we conclude
that
\begin{align*}
\ln r(\psi T) &=\frac{1}{p}\ell\bigl(\vert\psi\vert^p\bigr)= \frac{1}{p}\max_{\mu\in M_\delta (\cal
C)}\left(\int_{{X}}\ln\vert\psi\vert^p\,d\mu +{\tau(\mu)}\right)\\[6pt]
&= \max_{\mu\in M_\delta (\cal C)}\left(\int_{{X}}\ln \vert\psi\vert\,d\mu
 +\frac{\tau(\mu)}{p}\right). \qed
\end{align*}

\begin{rem}\label{13..7}
1) If in Definition \ref{13..3}\ \ $A\!:{\cal C}\to{\cal C}$ is an invertible conditional
expectation operator then in view  of Proposition \ref{8..3} the Variational Principle established
coincides with the variational principle (\ref{1,,4}).

2) If  $p\to \infty$ then  formula  (\ref{13,,15})  transforms into formula (\ref{1,,4}) and this
restores formally the variational principle (\ref{1,,4}) for the spaces of $C(X)$ and $L^\infty(X)$
type for an arbitrary $\delta$ and any transfer operator $A$ for $({\cal C},\delta)$ (cf.
Introduction).

3) To continue the previous remark we observe that if instead of $A$ mentioned in Definition
\ref{13..3} one takes a certain other transfer operator $A^0\!:\cal C\to \cal C$ such that
$A=A^0\varphi$ with some $0\le \varphi\in{\cal C}$ (for instance in the Model Example considered
above one can set $A^0=E$ and $\varphi = (d\beta(m)/dm)\circ\beta$ (see (\ref{13,,14}) and
(\ref{13,,13}))) then looking through the proof of the theorem we see that formula (\ref{13,,15})
transforms into
$$
 \ln r(\psi T) = \max_{\mu\in M_\alpha(\cal C)}\left(\int_{{X}}\biggl(\ln\vert\psi\vert +
 \frac{\ln\vert\varphi\vert}{p}\biggr)\, d\mu +\frac{\tau(\mu)}{p}\right).
$$
If $p\to\infty$ then this formula also transforms into formula (\ref{1,,4}). So weighted shift
operators in $C(X)$ and  $L^\infty(X)$ spaces `do not care' about transfer operators.

4) In the proof of Theorem \ref{13..1} we have used only condition b) of Definition \ref{13..3} and
have not exploited condition a) at all. We have inserted condition a) into Definition~\ref{13..3}
simply to emphasize the relation between the shift operator $T$ and endomorphism~$\delta$.

5) Recalling Remarks \ref{13..4} \, 3) \, and \ref{7..4} \, 4) \, we observe that even in the
standard space $L^1 (Y,m)$ an abstract weighted shift operator in general is \emph{not} generated
by any measurable mapping $\beta$ as in the Model Example considered above. So even in this
situation Theorem~\ref{13..6} is a generalization of the corresponding result form \cite{Bakh}.
\end{rem}

\end{document}